\documentclass[a4paper,12pt]{article}

\usepackage{amsfonts,amssymb,amsbsy,amsmath,amsthm,mathrsfs,enumerate,verbatim}

  \usepackage{graphics} 
  \usepackage{epsfig} 
  \usepackage{graphicx}
\usepackage{algorithmic,algorithm}
  \usepackage{epstopdf}

\usepackage{amsfonts}
\usepackage{amsmath}
\usepackage{amssymb}
\usepackage{graphicx}
\usepackage{epsf}
\usepackage{epsfig}

\usepackage{color}
\usepackage{afterpage}
\usepackage{verbatim}

\usepackage{times}

\newtheorem{theo}{Theorem}

\numberwithin{equation}{section}

\theoremstyle{definition}

\newcommand{\bel}{\begin{equation} \label}
\newcommand{\ee}{\end{equation}}

\def\beq{\begin{equation}}
\def\eeq{\end{equation}}
\newcommand{\bea}{\begin{eqnarray}}
\newcommand{\eea}{\end{eqnarray}}
\newcommand{\beas}{\begin{eqnarray*}}
\newcommand{\eeas}{\end{eqnarray*}}

{

 \usepackage{graphicx,color}
 \definecolor{mygreen}{cmyk}{1,0,1,0.1}

\usepackage{epsf}
\usepackage{epstopdf}

\usepackage{pdfsync}

\graphicspath{
  {Figures/}
  {/chalmers/users/larisa/PAPERS/BeilinaGajnova/ARXIV2019/PNG/Test1a/}
  {/chalmers/users/larisa/PAPERS/BeilinaGajnova/ARXIV2019/PNG/Test1b/}
  {/chalmers/users/larisa/PAPERS/BeilinaGajnova/ARXIV2019/PNG/Test2a/}
 {/chalmers/users/larisa/PAPERS/BeilinaGajnova/ARXIV2019/PNG/Test2b/}
}

\begin{document}

\title{Time-adaptive optimization in a parameter identification
  problem of HIV infection}

\author{L.~Beilina  \thanks{Department of Mathematical Sciences, Chalmers University of Technology and University of Gothenburg, SE-42196 Gothenburg, Sweden, e-mail:
\texttt{\
larisa@chalmers.se}}
\and
I.~Gainova \thanks{Sobolev Institute of Mathematics, 630090, Novosibirsk, Russia, e-mail: \texttt{\
gajnova@math.nsc.ru}}
}

\date{}

\maketitle

\begin{abstract}

The paper considers a time-adaptive method for determination of drug efficacy in a parameter identification problem (PIP) for system of ordinary differential equations (ODE) which describe dynamics of the primary HIV infection. Optimization approach to solve this problem is presented and a posteriori error estimates in the Tikhonov functional and Lagrangian are formulated. Based on these estimates a time adaptive   algorithm is formulated and numerically tested for different scenarios of noisy observations of virus population function. Numerical results show significant improvement of reconstruction of drug efficacy parameter when using time adaptive mesh refinement compared to usual gradient method applied on a uniform time mesh.

\end{abstract}

\maketitle


\section{Introduction}

\label{sec:0}

Parameter identification problems are frequently occurring within
biomedical applications. These problems are often non-linear and
ill-posed, and thus challenging to solve numerically.  For efficient
solution of parameter identification problems  a time-adaptive
method was recently proposed \cite{Sprg-1}. This method uses ideas
of an adaptive finite element method for solution of different
coefficient inverse problems for partial differential equations
(PDE) and has shown that it significantly improves reconstruction of
parameters \cite{BJ, BR, Beilina1, BJ2, BookBK}. The main idea of an
adaptive method is to  minimize a Tikhonov functional on  locally
refined meshes via a posteriori error estimates for the finite
element approximation of the solution of an inverse problem under
investigation.

This work is a continuation of the work by authors \cite{Sprg-2}
where they studied the model proposed in \cite{Sriv09} describing
the effect of the drug therapy on the dynamics of the Human
Immunodeficiency Virus (HIV) infection. In \cite{Sprg-2} a
time-adaptive method was formulated to determine the drug efficacy
in mathematical model of HIV infection using measurements in time of
all functions in this ODE system. Numerical simulations were not
presented in \cite{Sprg-2}. In the current work we consider a more
realistic case when only the virus population function is measured
and present numerical results of time-adaptive reconstruction of
drug efficacy from noisy measurements of virus population function
on an initial non-refined mesh. New a posteriori error estimate
between regularized and computed parameters is presented. Based on
this estimate, a time adaptive algorithm is formulated and
numerically tested on the reconstruction of drug efficacy from noisy
measurements of virus population function.

The time-adaptive method proposed in this paper can eventually be
used by clinicians to determine the drug-response for each treated
individual. Mathematical modelling helps to understand the biological mechanisms underlying in the base of action of antiviral drugs \cite{Boch2012}. The exact knowledge of the personal drug efficacy can aid in the determination of the most suitable drug as well as the
most optimal dose to an individual, in the long run resulting in a
personalized treatment with maximum efficacy and minimum adverse
drug reactions.

The outline of the paper is as follows.
The short biological description of the mathematical model
 is given in section \ref{sec:1-0}.
In section \ref{sec:1} the forward
and parameter identification problems
 are formulated.
The optimization method to solve the parameter identification  problem
is presented in  section \ref{sec:2}.
The finite element method
is formulated in section \ref{sec:fem}  and a posteriori error estimates
are presented in section \ref{sec:apostframework}.
An adaptive algorithm for solution of PIP  is formulated In section \ref{sec:fem_IP}. Finally,
in section \ref{sec:6}  numerical examples  confirm
the proposed time-adaptive algorithm.

\section{The mathematical model and its biological description}

\label{sec:1-0}

Despite the efforts by the international community to eradicate HIV
infection, the problem of its transmission, treatment and quality of
life of people living with HIV remains actual. According to
materials presented on the VI Eastern Europe and Central Asia AIDS
Conference and the latest data on HIV (UNAIDS, 2018), there are
currently more than 37 million people living with HIV globally and
an estimated two million new infections were recorded every year (http://aidsinfo.unaids.org).

In general case for living organisms the genetic information goes from the storage in DNA through messenger RNA (mRNA) to protein synthesis in the ribosomes. The process of converting the genetic information from DNA to mRNA is called \emph{transcription} \cite{Martin}. In the case of retroviruses, such as HIV, HIV's genetic information is encoded in form of RNA. HIV inserts its RNA into the host cell. Here viral RNA is reversely transcribed into HIV DNA, which is compatible with genetic material of the host cell (\emph{reverse transcription}). This DNA is transported to the cell's nucleus and incorporated into the DNA of the infected cell (\emph{integration}). To perform the reverse transcription of RNA into DNA, HIV carries its own enzyme called
\emph{reverse transcriptase}, that catalyzes the reverse transcription. Antiviral drugs inhibiting this enzyme (called \emph{Reverse Transcriptase Inhibitors})
will be able to prevent the production of new viruses  \cite{Cher2013, Patrick}.

Our basic mathematical model in this work is the model proposed in
\cite{Sriv09} which describes the effect of Reverse Transcriptase
Inhibitor (RTI) on the dynamics of HIV infection.  In this model the infected class
of CD4+ T-cells is subdivided into two subclasses: pre-RT class and
post-RT class. Pre-RT class consists of the infected CD4+ T-cells in
which reverse transcription is not completed, and post-RT class
consists of those infected CD4+ T-cells where the reverse
transcription is completed such that they are capable to produce
virus.
The mathematical model is:
\begin{equation}
  \label{state_forw}
\begin{array}{lcl}
\frac{du_1}{dt}   &=&{f_1(u(t), \eta(t))}=s -  k u_1(t) u_4(t) - \mu u_1(t) +(\eta(t) \alpha + b) u_2(t),\\[8pt]
\frac{du_2}{dt}   &=&{f_2(u(t), \eta(t))}=k u_1(t) u_4(t) - (\mu_1+ \alpha + b )u_2(t),\\[8pt]
\frac{du_3}{dt}    &=&{f_3(u(t), \eta(t))}=(1-\eta(t))\alpha u_2(t) -
\delta u_3(t),\\[8pt]
\frac{du_4}{dt}   &=&{f_4(u(t), \eta(t))}=N \delta u_3(t) - c u_4(t),
\end{array}
\end{equation}
with initial conditions
\begin{equation}\label{inidata}
\begin{array}{ll}
  u_1(0)=u^0_1=300\ mm^{-3},& u_2(0)=u_2^0=10\ mm^{-3}, \\[4pt]
  u_3(0)=u^0_3=10\ mm^{-3}, &u_4(0)=u^0_4=10\ mm^{-3}.
\end{array}
\end{equation}
Throughout the paper we denote by $\Omega_{T}= [0,T]$ the time domain
for $T>0$, where $T$ is the final observation time. In system
(\ref{state_forw}) the functions $u_i, i=1,2,3,4$ are defined as follows:
\begin{itemize}
\item $u_1(t)$ -- uninfected target cells population,
\item $u_2(t)$ -- infected target cells from pre-RT class,
\item $u_3(t)$ -- infected target cells from post-RT class,
\item $u_4(t)$ is the virus population function.
\end{itemize}
The initial data \eqref{inidata} are chosen
such that they satisfy two steady states (see details in
\cite{Sriv09}).

The system \eqref{state_forw} can be presented in the following compact form:
\begin{eqnarray}
\frac{du}{dt} &=& f(u(t),\eta(t)) ~~t \in [0,T], \label{forward_bio1} \\
u(0) &=& u^0, \label{forward_bio2}
\end{eqnarray}
with
\begin{equation}
  \begin{split}
    u=u(t) &= (u_1(t), u_2(t), u_3(t), u_4(t))^T, \\
    u^0 &= (u_1(0), u_2(0), u_3(0), u_4(0))^T, \\
    \frac{du}{dt} &= \left(\frac{\partial u_1}{ \partial t}, \frac{\partial u_2}{ \partial t},\frac{\partial u_3}{ \partial t}, \frac{\partial u_4}{ \partial t} \right)^T, \\
    f(u(t),\eta(t)) &= (f_1, f_2, f_3, f_4)(u(t),\eta(t))^T \\
    &= ( f_1(u_1,...,u_4, \eta(t)), ...,f_4(u_1, ..., u_4, \eta(t)))^T.
\end{split}
  \end{equation}

\newpage
\medskip
\textbf{Table 1}
\medskip

{\small
\begin{tabular}{llll}
  \hline
  Parameter & Value & Units & Description \\ \hline\\
  $s$ & $10$ & $mm^{-3}day^{-1}$ & inflow rate of T cells \\[4pt]
  $\mu$ & $0.01$ & $day^{-1}$ &natural death rate of T cells  \\[4pt]
  $k$ & $2.4$E-5  & $mm^{3}day^{-1}$ & interaction-infection rate of  T cells \\[4pt]
  $\mu_1$ & $0.015$ & $day^{-1}$ & death rate of infected  cells \\[4pt]
  $\alpha$ & $0.4$ & $day^{-1}$ & transition rate from pre-RT infected T cells class to post-RT class\\[4pt]
  $b$ & $0.05$ & $day^{-1}$ & reverting rate of infected cells return to uninfected class \\[4pt]
  $\delta$ & $0.26$ & $day^{-1}$ & death rate of actively infected cells \\[4pt]
  $c$ & $2.4$ & $day^{-1}$ & clearance rate of virus\\[4pt]
  $N$ &  $1000$ & ${vir}/{cell}$ & total number of viral particles produced by an infected cell \\
  &&&\\
  \hline
\end{tabular}
}

\section{The mathematical model and parameter identification problem}

\label{sec:1}

In the model \eqref{state_forw} we assume that $f \in C^1(\Omega_T)$
is Lipschitz continuous and  the function $\eta(t) \in
 C(\Omega_T)$ represents the unknown drug efficacy  which belongs to the set of
 admissible functions $M_{\eta}$:
\begin{equation}\label{6.1}
M_{\eta} = \{ \eta(t) : \eta(t)~ \in~ \left[0 , 1 \right]~
\textrm{in} ~\Omega_T ,~ \eta(t) =0 \textrm{ outside of } \Omega_T
\}.
\end{equation}

To formulate the parameter identification problem we assume that all
parameters in system (\ref{state_forw}) are known except the
parameter $\eta(t)$ which describes efficacy of the drug. The
typical values of  parameters $\{s, \mu, k, \mu_1, \alpha, b,
\delta, c, N \}$ in \eqref{state_forw} are  taken from
\cite{Sriv09}, see Table~1.

\medskip

\textbf{Parameter Identification Problem (}PIP\textbf{)}. Assume
that conditions (\ref{6.1}) hold and  parameters \linebreak $\{s $,
$\mu$, $k$, $\mu_1$, $\alpha$, $b$, $\delta$, $c$, $N \}$ in system
(\ref{state_forw})   are known. Assume further that the function $\eta(t) $
is unknown inside the domain $\Omega_T$. The PIP is: determine $\eta(t)$ for $t
\in \Omega_T,$  under the condition  that the  virus
population function $g(t)$ is known
\begin{equation}
u_4(t) = g(t),~~ t \in [T_1,T_2], 0 \leq T_1 < T_2 \leq T.  \label{6.4}
\end{equation}
Here, the function $g\left(t\right) $ presents observations of the function $u_4\left(t\right)$
 inside  the observation interval $[T_1,T_2]$.

 Note, that we solve the PIP on the time interval $[0,T]$ and assume
 that observations of $g(t)$ can even be on the more narrow interval
 $[T_1,T_2] \subset [0,T] $.  Numerical results of section \ref{sec:6}
 show that  reconstruction of the
 parameter $\eta(t)$ is not very good on the time interval where observations are not
 available  and thus,   observations of the virus population function $u_4(t)$
 should be taken as early as possible from  the date when the virus  started to be reproduced in the body of the host.

\section{Optimization method}

\label{sec:2}

Let $H$ be a Hilbert space of functions defined in $\Omega_T$.
To determine $\eta(t)$, $t\in [0,T]$ in PIP we minimize the  following Tikhonov functional
\begin{equation}\label{Tikh_f0}
J(\eta)=\frac{1}{2}   \int\limits_{T_1}^{T_2}(u_4(t)-g(t))^{2}z_{\zeta}\left(  t\right)
~\mathrm{d}t+\frac{1}{2}\gamma\int\limits_{0}^{T}(\eta-\eta^{0})^{2}dt.
\end{equation}
Here, the solution  $u_4(t)$ of the
system (\ref{state_forw})  with parameter $\eta(t)$, $g(t)$ is the observed virus population function, $\eta^0$ is
the initial guess for the parameter $\eta(t)$ and $\gamma\in (0,1)$ is
the regularization parameter, $z_{\zeta}(t), \zeta \in \left( 0,1\right) $ is smoothness
function which can be defined similarly  to \cite{Sprg-2}.

 To  find the function $\eta(t) \in H$ which minimizes the Tikhonov functional (\ref{Tikh_f0}) we
 seek for a stationary point of (\ref{Tikh_f0})  with respect to $\eta$ which satisfies
\begin{equation}\label{minimum}
J^{\prime }(\eta)(\bar{\eta})=0, ~~\forall \bar{\eta} \in H.
\end{equation}

To find minimum of (\ref{Tikh_f0}) we use the Lagrangian approach and
introduce the Lagrangian
\begin{equation}\label{Lagran_eta}
  L(v)=J(\eta)+  \sum_{i=1}^4 \int\limits_{0}^{T} \lambda_i
  \left(\frac{du_i}{dt} - f_i \right) ~dt,
\end{equation}
where $u(t)=(u_1(t),u_2(t),u_3(t),u_4(t))$ is the solution of the system (\ref{state_forw}), $\lambda(t)$ is the Lagrange multiplier $\lambda(t)=(\lambda_1(t),\lambda_2(t),\lambda_3(t),\lambda_4(t))$ and $v= (\lambda,u,\eta)$.

 Let us introduce  following spaces needed for further analysis
\begin{equation}
\begin{array}{rl}
H_{u}^{1}(\Omega_T)  &  =\{f\in H^{1}(\Omega_T):f(0)=0\},\\
H_{\lambda}^{1}(\Omega_T)  &  =\{f\in H^{1}(\Omega_T):f(T)=0\},\\
U  &=H_{u}^{1}(\Omega_T)\times H_{\lambda}^{1}(\Omega_T)\times C(\Omega_T),
\end{array}
\label{spaces}
\end{equation}
where all functions are real valued.

 To derive the Fr\'{e}chet derivative of the Lagrangian
(\ref{Lagran_eta})
we assume that functions $v=(\lambda,u,\eta)$ can be varied independently
of each other in the sense that
\begin{equation}\label{Frechet_der}
  L'(v)(\bar{v})=0,\quad \forall  \bar{v} =(\bar{\lambda},\bar{u},\bar{\eta}) \in U.
\end{equation}
Thus,  we consider  $L(v + \bar{v}) - L(v)$, single out the linear part with respect to $v$ of the obtained expression and neglect all nonlinear terms.
The   optimality condition \eqref{Frechet_der} means that for all
  $\bar{v} \in U$ we have
\begin{equation}
  L'(v; \bar{v}) = \frac{\partial L}{\partial \lambda}(v)(\bar{\lambda}) +  \frac{\partial L}{\partial u}(v)(\bar{u})
  + \frac{\partial L}{\partial \eta}(v)(\bar{\eta}) = 0,  \label{scalar_lagrang1}
\end{equation}
i.e., every component of \eqref{scalar_lagrang1} should be zero
out. Thus, the optimality conditions \eqref{Frechet_der}  yields
\begin{equation}\label{forward1}
\begin{split}
  0 = \frac{\partial L}{\partial \lambda}(v)(\bar{\lambda}) &=
- \alpha\int\limits_{0}^{T} u_2 (\lambda_1-\lambda_3) \bar{\eta} dt \\
 &+  \int\limits_{0}^{T} ( \dot{u_1}  - s + ku_1u_4 + \mu u_1 - (\eta\alpha + b) u_2) \bar{\lambda_1} dt \\
&+ \int\limits_{0}^{T} (\dot{u_2} - ku_1u_4 +  (\mu_1+ \alpha + b )u_2) \bar{\lambda_2} dt \\
&+\int\limits_{0}^{T}  (\dot{u_3} - (1-\eta)\alpha u_2 + \delta u_3) \bar{\lambda_3} dt \\
&+ \int\limits_{0}^{T} ( \dot{u_4} - N \delta u_3 + c u_4) \bar{\lambda_4} dt
 ~~~\forall \bar{\lambda} \in H_u^1(\Omega_T),
\end{split}
\end{equation}

\begin{equation} \label{control1}
\begin{split}
  0 = \frac{\partial L}{\partial u}(v)(\bar{u}) &=
-\int\limits_{0}^{T}(  \dot{\lambda}_{1} - \lambda_1 ku_{4} - \lambda_1\mu  + \lambda_2k u_{4}) \bar{u_1} dt  \\
&- \int\limits_{0}^{T}( \dot{\lambda}_{2} - \lambda_2(\mu_1+\alpha+b) + \lambda_1 (\eta \alpha + b) + (1 - \eta)\alpha \lambda_3)  \bar{u_2} dt
 \\
 &-\int\limits_{0}^{T}(   \dot{\lambda}_{3} - \lambda_3\delta + \lambda_4N\delta ) \bar{u_3} dt \\
&- \int\limits_{0}^{T}(  \dot{\lambda}_{4} - \lambda_4 c - \lambda_1k u_1 + \lambda_2ku_1 ) \bar{u_4} dt
+ \int\limits_{T_1}^{T_2} (u_4 - g)z_{\zeta} \bar{u_4} dt
~~~\forall \bar{u} \in  H_{\lambda}^1(\Omega_T),
\end{split}
\end{equation}
\begin{equation} \label{grad1new}
\begin{split}
0 &= \frac{\partial L}{\partial \eta}(v)(\bar{\eta})
= \gamma\int\limits_{0}^{T}(\eta-\eta^{0})\bar{\eta} dt
+ \alpha \int_0^T u_2(\lambda_3 - \lambda_1)\bar{\eta} dt
~~~\forall \bar{\eta} \in C\left(\Omega_T\right).
\end{split}
\end{equation}

\medskip

The  equation  \eqref{forward1} corresponds to the forward problem (\ref{state_forw})-(\ref{inidata}), the  equation \eqref{control1} --- to  the  following adjoint problem
\begin{equation}\label{adjoint}
\begin{array}{lcl}
  \frac{\partial \lambda_1}{\partial t} & = & \lambda_1 ku_{4} + \lambda_1\mu  - \lambda_2k u_{4}, \\[8pt]
  \frac{\partial \lambda_2}{\partial t} & = & \lambda_2(\mu_1+\alpha+b) - \lambda_1 (\eta \alpha + b)  - (1-\eta) \alpha \lambda_3, \\[8pt]
  \frac{\partial \lambda_3}{\partial t} & = & \lambda_3\delta - \lambda_4N\delta, \\[8pt]
   \frac{\partial \lambda_4}{\partial t} & = &  \lambda_4 c + \lambda_1k u_1 - \lambda_2ku_1 + (u_4 - g)z_{\zeta} , \\[4pt]
\lambda_i(T) &=&0,\quad i=1,\ldots,4.
\end{array}
\end{equation}
which can be rewritten in  the  compact form as
\begin{equation}\label{adjoint}
  \begin{split}
    \frac{\partial \lambda}{\partial t} &=  \tilde{f}(\lambda(t)),\\
    \lambda_i(T) &= 0,\quad i=1,\ldots,4
    \end{split}
\end{equation}
with
\begin{equation}
  \begin{split}
    \lambda= \lambda(t) &=(\lambda_1(t), \lambda_2(t), \lambda_3(t), \lambda_4(t))^T, \\
    0&= (\lambda_1(T), \lambda_2(T), \lambda_3(T), \lambda_4(T))^T, \\
    \frac{d\lambda}{dt} &=  \left(\frac{\partial \lambda_1}{ \partial t}, \frac{\partial \lambda_2}{ \partial t},\frac{\partial \lambda_3}{ \partial t}, \frac{\partial \lambda_4}{ \partial t} \right)^T, \\
    \tilde{f}(\lambda(t)) &=( \tilde{f}_1,  \tilde{f}_2,  \tilde{f}_3,  \tilde{f}_4)( \lambda(t))^T.
\end{split}
  \end{equation}
The adjoint system should be solved backwards in time with already known solution $u(t)$ to the forward problem (\ref{state_forw})-(\ref{inidata}) and a given  measurement function $g(t)$.

For the case when $u$ and $\lambda$ are exact solutions of the forward
(\ref{state_forw})-(\ref{inidata}) and adjoint (\ref{adjoint})
problems, respectively, to the known function $\eta$,
we get from \eqref{Lagran_eta} that
\begin{equation}
L(v(\eta)) = J(\eta),
\end{equation}
and thus the Fr\'{e}chet derivative of the Tikhonov functional can be
written as
\begin{equation}\label{derfunc}
\begin{split}
 J'(\eta) := &J_\eta(u(\eta), \eta) =  \frac{\partial J}{\partial \eta}(u(\eta), \eta)
  =  \frac{\partial L}{\partial \eta}(v(\eta)).
\end{split}
\end{equation}
Using (\ref{grad1new})  in  (\ref{derfunc}),
we get the following expression for the Fr\'{e}chet derivative of the Tikhonov functional
\begin{equation}\label{Frder2}
  J'(\eta)(t)  =\gamma  (\eta-\eta^{0})(t) + \alpha  u_2 (\lambda_3 -
   \lambda_1)(t) = 0,
\end{equation}

Thus, to find  the unknown parameter $\eta$ which minimizes the Tikhonov functional (\ref{Tikh_f0})  we can use the following expression
\begin{equation} \label{coef2}
\begin{split}
\eta &= \frac{1}{\gamma}
 \alpha u_2 (\lambda_1-\lambda_3) + \eta^{0}.
\end{split}
\end{equation}

\section{ Finite Element Discretization}
\label{sec:fem}

For  solution of \eqref{Frechet_der} we will use the finite element discretization and    consider a partition $\mathcal J_{\tau} = \{J\}$ of the time
domain $\Omega_T = [0,T]$
 into time subintervals $J=(t_{k-1},t_k]$
  of the  time step $\tau_k =  t_k - t_{k-1}$. We  define also the piecewise-constant time-mesh
  function $\tau$ such that
\begin{equation}\label{neshfunction}
\tau(t) = \tau_k, ~~\forall J \in  J_{\tau}.
\end{equation}

For discretization of the state and adjoint problems we   define the
finite element spaces $W_{\tau}^{u}\subset H_{u}^{1}\left(\Omega_T\right)
$ and $W_{\tau}^{\lambda}$ $\subset H_{\lambda} ^{1}\left(\Omega_T\right)
$ for $u$ and $\lambda$, respectively, as
\begin{equation}
\begin{array}{rl}
W_{\tau}^{u}  &=\{f\in  H_{u}^{1}:  f|_J \in P^1(J)~~~ \forall J \in J_{\tau}\},\\
W_{\tau}^{\lambda} &=\{f\in  H_{\lambda}^{1}: f|_J \in P^1(J)~~~ \forall J \in J_{\tau} \}.
\end{array}
\label{femspaces}
\end{equation}

For the function $\eta(t)$ we also introduce the finite element space $W_{\tau}^{\eta}\subset L_2\left(\Omega_T\right) $ consisting of piecewise
constant functions
\begin{equation}
W_{\tau}^{\eta} =\{f\in  L_{2}\left( \Omega_T\right) :  f|_J \in P^0(J)~~ \forall J \in J_{\tau}\}.
\end{equation}

We use different finite element spaces since we are working in a finite
dimensional space and all norms in finite dimensional spaces are
equivalent. Next we denote $U_{\tau}= W_{\tau}^{u}\times W_{\tau}^{\lambda}\times W_{\tau}^{\eta}$
such that $U_{\tau}\subset U$.

Now the finite element method for (\ref{Frechet_der}) is: find $v_{\tau}\in U_{\tau}$ such that
\begin{equation}
L^{\prime}\left(  v_{\tau};\bar{v}\right)  =0,~~\forall\overline{v}\in U_{\tau}.
\label{discr_lagr_ac}
\end{equation}

Since the forward \eqref{state_forw} - \eqref{inidata} and adjoint \eqref{control1}  problems are nonlinear their solutions
can be found by Newton's method.
For the discretization
\begin{equation*}
\frac{\partial u}{\partial t} = \frac{u^{k+1} - u^k}{\tau_k}
  \end{equation*}
the variational formulation of the forward problem \eqref{state_forw} - \eqref{inidata}  for all $\bar{u} \in H_u^1(\Omega_T)$  is:
\begin{equation}\label{4.26}
(u^{k+1}, \bar{u}) -(u^k,  \bar{u}) - \tau_k f(u^{k+1},  \bar{u}) = 0.
  \end{equation}
Denoting
\begin{equation}\label{4.26a}
  \begin{split}
  \tilde{u} &= u^{k+1}, \\
  F(\tilde{u}) &= \tilde{u} - \tau_k f(\tilde{u}) - u^k
\end{split}
  \end{equation}
  we can rewrite \eqref{4.26} as
\begin{equation}\label{4.27}
(F(\tilde{u}),  \bar{u}) = 0.
  \end{equation}
For solution $F(\tilde{u})=0$  the Newton's method can be used for the iterations $n=1,2,...$  \cite{Burden}
\begin{equation}\label{4.31}
\tilde{u}^{n+1} = \tilde{u}^n - [ F'(\tilde{u}^n)]^{-1} \cdot F(\tilde{u}^n).
  \end{equation}
Here,  we can determine  $F'(\tilde{u}^n)$ via  definition of  $F(\tilde{u})$  in \eqref{4.26a} as
\begin{equation*}
F'(\tilde{u}^n) = I - \tau_k f'(\tilde{u}^n),
 \end{equation*}
where $I$ is the identity matrix, $f'(\tilde{u}^n)$ is the Jacobian of
$f$  (the right hand side of the forward problem   \eqref{state_forw}) at $\tilde{u}^n$ and $n$ is the iteration number in Newton's method. We note that the finite element method \eqref{discr_lagr_ac}
will work even in this case, see   details in  \cite{Johnson}.

In a similar  way the Newtons's method  can be derived for the solution the adjoint problem \eqref{adjoint}.  Since we solve the adjoint problem backwards in time, we  discretize
time derivative as
\begin{equation}\label{4.32}
\frac{\partial \lambda}{\partial t} = \frac{\lambda^{k+1} - \lambda^{k}}{\tau_k}
  \end{equation}
for the already known  $\lambda^{k+1}$ values, and write the variational formulation of the adjoint problem for all
$\bar{\lambda} \in H_\lambda^1(\Omega_T)$
\begin{equation}\label{4.33}
( \lambda^{k}  - \lambda^{k+1} + \tau_k \tilde{f}(\lambda^k) ,\bar{\lambda}) = 0.
  \end{equation}
Denoting
\begin{equation}\label{4.34}
  \begin{split}
  \tilde{\lambda} &= \lambda^{k}, \\
  \tilde{F}(\tilde{\lambda}) &= \tilde{\lambda} + \tau_k \tilde{f}(\tilde{\lambda})  -  \lambda^{k+1},
\end{split}
  \end{equation}
  we can rewrite \eqref{4.33}  for all
$\bar{\lambda} \in H_\lambda^1(\Omega_T)$ as
\begin{equation}\label{4.35}
(\tilde{F}(\tilde{\lambda}),  \bar{\lambda}) = 0.
  \end{equation}
For solution $ \tilde{F}(\tilde{\lambda})=0$  we use again Newton's method  for iterations $n=1,2,...$
\begin{equation}\label{4.36}
\tilde{\lambda}^{n+1} = \tilde{\lambda}^n - [ \tilde{F}'(\tilde{\lambda}^n)]^{-1} \cdot \tilde{F}(\tilde{\lambda}^n).
  \end{equation}
We compute  $\tilde{F}'(\tilde{\lambda}^n)$ using the  definition of  $\tilde{F}(\tilde{\lambda})$  in \eqref{4.34} as
\begin{equation*}
 \tilde{F}'(\tilde{\lambda}^n) = I + \tau_k \tilde{f}'(\tilde{\lambda}^n),
 \end{equation*}
where $I$ is the identity matrix, $\tilde{f}'(\tilde{\lambda}^n)$ is the Jacobian of
$\tilde{f}$ (the right hand side of the adjoint problem \eqref{adjoint}) at $\tilde{\lambda}^n$ and $n$ is the iteration number in Newton's method.

\section{A Posteriori Error Estimates }

\label{sec:apostframework}

We consider the function $\eta \in C(\Omega_T)$ as a minimizer of the
Lagrangian \eqref{Lagran_eta}, and $\eta_{\tau} \in W_{\tau}^\eta$ its  finite element
approximation.
Let us assume that we know good approximation to the exact solution
$\eta^*  \in  C(\Omega_T)$.  Let $g^{\ast }(t)$ be the exact data and the function
$g_{\sigma }(t)$ represents the error level in these data. We assume
that measurements $g(t)$ in (\ref{6.4}) are given with some noise
level (small) $\sigma$ such that
\begin{equation}
g(t)=g^{\ast }(t)+g_{\sigma}(t);\,~ g^{\ast },g_{\sigma }\in L_{2}\left( \Omega_{T}\right) ,~ \left\Vert g_{\sigma }\right\Vert _{L_{2}\left( \Omega_{T}\right) }\leq \sigma .
\label{6.18}
\end{equation}

Accordingly \cite{BKK} we assume that
\begin{equation}\label{6.19}
\gamma = \gamma(\sigma) = \sigma^{2 \mu},~~\mu \in(0,1/4), ~~\sigma \in (0,1)
\end{equation}
and
\begin{equation}\label{6.20}
\| \eta_0 - \eta^* \| \leq \frac{\sigma^{3\mu}}{3},
\end{equation}
where $\eta^*$ is the exact solution of PIP with the exact data $g^*(t)$. Let
\begin{equation}
V_\varepsilon(\eta) = \{ x \in  C(\Omega_T): \| \eta - x \| < \varepsilon ~~~\forall \eta \in C(\Omega_T)  \}.
  \end{equation}
Assume  that
for all $\eta \in V_1(\eta^*)$ the operator
\begin{equation}
F(\eta) =   \frac{1}{2}   \int\limits_{T_1}^{T_2}(u_4(\eta, t)-g(t))^{2}z_{\zeta}\left(  t\right)~\mathrm{d}t
 \end{equation}
has the
Fr\'{e}chet derivative  $F'(\eta)$ which is bounded and Lipshitz continuous
in  $V_1(\eta^*)$ for $D_1, D_2 = const. > 0$
\begin{equation}\label{6.21}
  \begin{split}
  \| F'(\eta) \| &\leq D_1 ~ ~~ \forall \eta \in V_1(\eta^*), \\
  \| F'(\eta_1) -  F'(\eta_2)  \| &\leq D_2 \| \eta_1 - \eta_2 \|  ~~ \forall \eta_1, \eta_2 \in  V_1(\eta^*).
\end{split}
  \end{equation}

\subsection{An a posteriori error estimate for the Tikhonov functional}

\label{sec:errorfunc}

In the Theorem 1 we derive an a posteriori error estimate for the error in
the Tikhonov functional (\ref{Tikh_f0}) on the finite element time partition $\mathcal J_{\tau}$.

\textbf{Theorem 1}. \emph{\ We assume that there exists minimizer $\eta \in
  C(\Omega_T)$ of the functional $J(\eta)$ defined by (\ref{Tikh_f0}). We assume also that
  there exists finite element approximation of a minimizer $\eta_{\tau} \in W_{\tau}^{\eta}$ of
  $J(\eta)$.  Then the following
  approximate a posteriori error estimate for the error $ e=|| J(\eta) - J(\eta_{\tau}) ||_{L^2(\Omega_T)}$ in the Tikhonov functional (\ref{Tikh_f0}) holds  true}
\begin{equation}\label{theorem2}
  e= || J(\eta) - J(\eta_{\tau}) ||_{L^2(\Omega_T)} \leq C_I C \left\| J^{\prime }(\eta_{\tau})\right\| _{L^2(\Omega_T)}  || \tau \eta_{\tau} ||_{L_2(\Omega_T)}
\end{equation}
with positive constants $C_I, C > 0$ and where
\begin{equation}
J^{\prime }(\eta_{\tau}) = \gamma (\eta_{\tau} -\eta^{0}) - \alpha {u_2}_{\tau} ({\lambda_1}_{\tau}
- {\lambda_3}_{\tau}). \label{theorem2_1}
\end{equation}

\textbf{Proof}

Proof follows from the Theorem 5 of \cite{Sprg-1}.

$\square$

\subsection{A posteriori error estimate of the minimizer on  refined meshes}

\label{sec:adaptrelax}

Theorems 2  and 3  present two a posteriori  error estimates for a minimizer $\eta$
of  the functional \eqref{Tikh_f0}.

\textbf{Theorem 2}

 \emph{ Let $\eta_\tau \in W_\tau^\eta$
   be a finite element approximation on the finite element mesh $J_\tau$ of the minimizer
   $\eta \in L^2(\Omega_T)$ of  the functional \eqref{Tikh_f0}
   with the mesh function $\tau(t)$. Then there exists a
  Lipschitz constant   $D = const.>0$ defined by
\begin{equation}
  \left\| J^{\prime }(\eta_1) - J^{\prime }(\eta_2) \right\|
  \leq D \left\| \eta_1 - \eta_2\right\| ,\forall \eta_1, \eta_2 \in L_2(\Omega_T),
  \label{2.10}
\end{equation}
  and interpolation constant $C_I$  independent on $\tau$
   such that the following
  a posteriori error estimate for the minimizer   $\eta $ holds  true
\begin{equation}
|| \eta_\tau  - \eta ||_{L_2(\Omega_T)} \leq  \frac{D}{\gamma } C_I  || \tau \eta_\tau ||_{L_2(\Omega_T)} ~ \forall \eta_\tau \in W_\tau^\eta.  \label{theorem1}
\end{equation}
 }

\textbf{Proof.}

 Proof  follows from the  Theorem 5.1 of \cite{KB}.

$\square$

\textbf{Theorem 3}

 \emph{Let $\eta_\tau \in W_\tau^\eta$
   be a finite element approximation on the finite element mesh $J_\tau$ of the minimizer
   $\eta \in L^2(\Omega_T)$ of  the functional \eqref{Tikh_f0}
   with the mesh function $\tau(t)$. Then there exists
  an interpolation constant $C_I$  independent on $\tau$
   such that the following
  a posteriori error estimate for the minimizer   $\eta $ holds}
\begin{equation}
|| \eta_\tau  - \eta ||_{L_2(\Omega_T)} \leq \sqrt{ \frac{\|R(\eta_\tau)\|}{\gamma } C_I  || \tau \eta_\tau ||_{L_2(\Omega_T)}} ~ \forall \eta_\tau \in W_\tau^\eta,  \label{theorem3}
\end{equation}
\emph{ where $R(\eta_\tau)$ is the residual defined as}
\begin{equation}\label{resid}
  R(\eta_\tau)(t)  = \gamma (\eta_\tau - \eta^0)(t) + \alpha  {u_2}_\tau
  ({\lambda_3}_\tau - {\lambda_1}_\tau)(t).
\end{equation}

\textbf{Proof.}

Proof follows from Theorems 1 and 2.

 $\square $

\section{Algorithms for solution of PIP}

\label{sec:fem_IP}

Here we present two algorithms for solution of PIP:
\begin{itemize}
  \item CGA -  usual  conjugate gradient algorithm (CGA) on a coarse time partition,
 \item ACGA - time-adaptive conjugate gradient algorithm which
   minimized the Tikhonov functional \eqref{Tikh_f0} on a locally
   refined meshes in time.
\end{itemize}

 We denote the nodal value of the
gradient at the observation points $\{t_i\}$ by $G^{m}(t_i)$ and compute
it accordingly to \eqref{Frder2} as
\begin{equation}
  G^{m}(t_i) = \gamma (\eta_\tau^m(t_i) - \eta_\tau^0(t_i)) +
  \alpha  {u_2}_\tau^m(t_i) ({\lambda_3}_\tau^m(t_i) - {\lambda_1}_\tau^m(t_i)).
\label{gradient2}
\end{equation}
The approximate computed solutions ${u_2}_\tau^{m} $and ${\lambda_{1,3}}_\tau^m$ are
obtained computationally by  Newton's method with
$\eta:={\eta_\tau}^{m}$.
A sequence $\{{\eta_\tau}^{m} \}_{m=1,...,M} $ of
approximations to $\eta$ is computed as  follows
\begin{equation}\label{conjgrad}
\begin{split}
\eta_\tau^{m+1}(t_i) &=  \eta_\tau^{m}(t_i)  + r^m d^m(t_i),
\end{split}
\end{equation}
with
\begin{equation*}
\begin{split}
 d^m(t_i)&=  -G^m(t_i)  + \beta^m  d^{m-1}(t_i),
\end{split}
\end{equation*}
and
\begin{equation*}
\begin{split}
 \beta^m &= \frac{|| G^m(t_i)||^2}{|| G^{m-1}(t_i)||^2},
\end{split}
\end{equation*}
where $d^0(t_i)= -G^0(t_i)$  and $G^{m}(t_i)$ is the gradient  vector  which is computed by (\ref{gradient2}) in time moments $t_i$.
In \eqref{conjgrad} the parameter $r^m$ is the step-size in the gradient update   at the iteration $m$ which is computed as
\begin{equation}
r^m = -\frac{(G^m, d^m)}{\gamma \| d^m\|^2}.
\end{equation}



\label{sec:adalg}

\begin{algorithm}[hbt!]
  \centering
  \caption{Conjugade Gradient  Algorithm  (CGA).\label{alg:cga}}
  \begin{algorithmic}[1]
    \STATE Choose time partition $\mathcal J_{\tau}$ of the time
  interval $(0, T)$. Start with the initial approximations
  ${\eta_\tau}^{0}$   and compute the sequence of ${\eta_\tau}^{m}$ for all $m >0$ in the following steps.
  \STATE Compute solutions $u_\tau^m = u_\tau\left(t,{\eta_\tau}^{m}\right)
  , \lambda_\tau^m = \lambda_\tau\left(t, {\eta_\tau}^{m}\right) $ of the state (\ref{state_forw})
  and adjoint (\ref{adjoint}) problems via  \eqref{4.31}, \eqref{4.36}, respectively, using Newton's method on the time partition
  $\mathcal J_{\tau}$.
 \STATE Compute gradient $G^{m}(t_i)$ on the time partition
  $I_{\tau}$  by (\ref{gradient2}).
\STATE Update the unknown parameter $\eta :=
  \eta_\tau^{m+1}$ using (\ref{conjgrad}) as
\begin{equation*}
\begin{split}
\eta_\tau^{m+1}(t_i) &=  \eta_\tau^{m}(t_i)  + r^m d^m(t_i).
\end{split}
\end{equation*}
  \STATE Compute residual  $R(\eta_\tau^m)$ using  \eqref{resid} with  solutions $u_\tau\left(t,{\eta_\tau}^{m}\right)
  ,\lambda_\tau\left(t, {\eta_\tau}^{m}\right) $ of the state (\ref{state_forw})
  and adjoint (\ref{adjoint}) problems.
\STATE For the tolerance  $0 <\theta < 1$ chosen by the user, stop computing the functions
  $\eta_\tau^{m}$ if either
  $||R(\eta_\tau^m)||_{L_2(\Omega_T)} \leq \theta$, or  norms of  residuals
  $||R(\eta_\tau^m)||_{L_2(\Omega_T)} $ abruptly grow, or norms of computed
  $||\eta_\tau^{m}||_{L_2(\Omega_T)} $ are stabilized. Otherwise, set $m:=m+1$ and go to Step 2.
  \end{algorithmic}
\end{algorithm}


\label{sec:7.1}

In the adaptive algorithm ACGA we have used  Theorem 3 for the error $e = \| \eta_\tau - \eta \|_{L_2(\Omega_T)}$ on locally refined meshes. More precisely,  first we  choose tolerance  $0 <\theta < 1$  and run adaptive algorithm until
\begin{equation*}
e = \| \eta_\tau - \eta \|_{L_2(\Omega_T)} \leq \theta.
\end{equation*}

For the time-mesh refinements we propose following  refinement procedure
based on the  Theorem 3.

\textbf{The} \textbf{Time Mesh Refinements Criterion}

\emph{Refine the time-mesh $\mathcal J_{\tau}$ in neighborhoods of those time-mesh
  points }$t\in {\Omega_T}$\emph{\ where the residual  }$\left\vert
R\left( \eta_{\tau}\right) \left( t\right) \right\vert
$\emph{\   defined in \eqref{resid} attains its maximal values. More precisely, let }$\beta _{1}\in \left(
0,1\right) $\emph{\ be the tolerance number. Refine the time-mesh in
  such subdomains of }${\Omega_T} $\emph{\ where}
\begin{equation*}
\left\vert R(\eta_\tau) \left( t\right)
\right\vert \geq \beta_{1}\max_{\Omega_T}\left\vert
 R(\eta_\tau) \left( t\right) \right\vert .
\end{equation*}

Using the above mesh refinement recommendation we propose the following time-adaptive algorithm in computations:

\begin{algorithm}[hbt!]
  \centering
  \caption{Adaptive Conjugate Gradient Algorithm (ACGA) .\label{alg:acga}}
  \begin{algorithmic}[1]
    \STATE  Peform steps 1-6 in  CGA algorithm. Let $M$ be the final iteration in CGA algorithm.

\STATE  Refine the time mesh $\mathcal J_{\tau}$ at all points where
\begin{equation}
\left\vert R(\eta_\tau^M) \left( t\right)
\right\vert \geq \beta_{1}\max_{\Omega_T}\left\vert
 R(\eta_\tau^M) \left( t\right) \right\vert .  \label{6.70}
\end{equation}
Here the tolerance number $\beta _{1}\in \left( 0,1\right) $ is chosen by
the user.
\STATE  Construct a new  time
partition $\mathcal J_{\tau}$ of the time interval $\left( 0,T\right) $.
 Interpolate the initial approximation $\eta_0$ from the
 previous  time partition to the new  time partition. Next,
 peform steps 1-6 in  CGA algorithm   on the new  time partition.
\STATE  Stop time partition refinements if norms  of residuals
$||R(\eta_\tau^M)||_{L_2(\Omega_T)} $  either increase or stabilize, compared with the
  previous time partition.
  \end{algorithmic}
\end{algorithm}

\section{Numerical results}
 \label{sec:6}

 In this section we present several numerical results  which show performance  and effectiveness of the time-adaptive reconstruction of  unknown parameter $\eta(t), t \in [0,T]$ in PIP  using ACGA algorithm.
Numerical results of reconstruction of function $\eta(t)$ using  usual conjugate gradient
 Algorithm 1 on the nonrefined time-meshes are presented in
 \cite{Martin}. We note that observations of all $u_i, i=1,2,3,4$
 functions in system \eqref{state_forw} were used in \cite{Martin} .

The goal of numerical tests of this note is to determine the unknown function $\eta(t)$
from observation of the virus
population function $u_4(t)$ in \eqref{state_forw} on the interval $[T_1,T_2] \subset  [0,T], 0 \leq T_1  < T_2 \leq T$.
 In all numerical tests  assumed that parameter $\eta(t)$ satisfy conditions
 \eqref{6.1} and is unknown in the system \eqref{state_forw}, but all
 other parameters $\{s $, $\mu$, $k$, $\mu_1$, $\alpha$, $b$,
 $\delta$, $c$, $N \}$ of this system are known and their values are
 chosen as in the Table 1.  The observation interval $[T_1, T_2]$ is such that  $T_2 = T = 300$, but $T_1$ is taken
differently in different tests since observations of the virus population function $u_4(t)$ can be
taken after the first $3-9$ weeks since the virus started to be
reproduced in the body of host.

For generation of data $u_4(t) = g(t)$  the
problem \eqref{state_forw}-\eqref{inidata}  was solved numerically with  exact values of the test
model function  $\eta(t)$. For solution of  problem \eqref{state_forw}-\eqref{inidata} was used Newton's method presented in section 4.

Next, the random  noise  was added  to the observed  solution $u_4(t)$ as
\begin{equation}
{u_4}_\sigma(t) = {u_4}_\sigma(t)( 1 + \sigma \alpha),
  \end{equation}
where $\sigma \in [0,1]$ is nose level and $\alpha \in [-1,1]$ is
random number.

In Algorithms  1, 2 it is of vital importance to take initial guess $\eta^0$
such that it satisfy condition  \eqref{6.20} which means that $\eta^0$ is located in the close
neighborhood of the exact solution. This condition is fullfilled in our PIP since we can compute explicitly values of the parameter $\eta(t)$ on the initial non-refined time mesh using, for example,  the  third  equation of system \eqref{state_forw} as
\begin{equation}
\eta(t) = 1 - \frac{\frac{\partial u_3(t)}{\partial t} + \delta u_3(t)}{\alpha u_2(t)}.
\end{equation}

We used following discretised version of this equation to get initial guess  $\eta_\tau^0$
\begin{equation}\label{approxeta}
\eta_\tau^0(t) \approx  1 - \frac{\frac{ {u_3}_\tau^{k+1} - {u_3}_{\tau}^{k}}{ \tau_k } + \delta {u_3}_\tau^k}{\alpha {u_2}_\tau^k}.
\end{equation}
Here, ${u_3}_\tau^{k+1}, {u_3}_\tau^k, {u_2}_\tau^k$ are known
computed approximations of functions $u_3, u_2$ at time iterations
$k+1$ and $k$, respectively.  We note that denominator is not
approaching zero because $\alpha = 0.4$ and ${u_2}_\tau(t) > 0 ~~ \forall t \in [0,T]$.
To get reasonable approximation  $\eta_\tau^0$ for the
 initial guess $\eta^0$ in Algorithm 2 we assume that noisy functions
${u_3}_{\sigma}, {u_2}_{\sigma}$ are known on the initial non-refined
 mesh, apply \eqref{approxeta}  and then use polynomial fitting to obtained noisy data $\eta_\tau^0$. Finally,  the condition \eqref{6.1}  was applied for the computed $\eta_\tau^0$  in order
to ensure that $\eta^0(t)$ belongs to the set of admissible
parameters $M_\eta$.

All tests are performed with tolerance $\Theta = 10^{-7}$ in ACGA algorithm and $\beta_1
= 0.1$ in \eqref{6.70}. The value of $\beta_1$ is chosen such that it
allows local refinements and avoids refinement of the very large time region in the time mesh.
 All tests are performed for  different $T_1= 25, 50, 100$  for the time interval $[T_1, T_2] = [T_1, 300]$ which corresponds to the fact that HIV virus can be detected in the first 3-9 weeks after infection.

Relative errors  in the reconstructed parameters $\eta(t)$ presented in the Tables  are   measured in
 $L_2$-norm and are computed as
 \begin{equation}\label{relerror}
e_{\eta} = \frac{\| \eta - \eta_\tau\|_{L_2(\Omega_T)}}{\| \eta \|_{L_2(\Omega_T)}}.
 \end{equation}

\subsection{Test 1}

\begin{table}[h!]
\center
\begin{tabular}{ |l|l|l|l|l| }
  \multicolumn{5}{c}{  $T_1 = 25$ }\\
  \hline
$\sigma$ & 5 \%  &   10\%   &  20\%   & 40\%       \\
nr.of ref. &   &    &   &        \\
\hline
$0$ &   0.1893  &  0.2022   &   0.2129 &    0.2203                \\
$1$  & 0.1151   &  0.1223 &   0.1279 &  0.2008 \\
$2$  &  0.0470   &  0.0391 &  & \\
$3$  &   0.0354  &  &  & \\
$4$  &  0.0242   &  &  & \\
\hline
\multicolumn{5}{c}{  $T_1 = 50$ }\\
  \hline
$\sigma$ & 5 \%  &   10\%   &  20\%   & 40\%       \\
nr.of ref. &   &    &   &        \\
\hline
$0$ &    0.1917 &   0.1933  &  0.1639   &    0.3498                \\
$1$  &  0.1194   &  0.1267 &  0.1027  &  0.2990  \\
$2$  &   0.0684  &  0.0550  &  0.1002  & 0.1755 \\
$3$  &   0.0337 &  0.0394 & 0.0657  &  0.1677  \\
$4$  &  0.0217 &  &  & \\
\hline
\multicolumn{5}{c}{  $T_1 = 100$ }\\
  \hline
$\sigma$ & 5 \%  &   10\%   &  20\%   & 40\%       \\
nr.of ref. &   &    &   &        \\
\hline
$0$ &   0.1560  &    0.1851    &   0.2494 &   0.3229                \\
$1$  &   0.1106  &   0.1275 & 0.1442   &   0.2035  \\
$2$  &   0.0775  &   0.0810  &   0.1132   &   0.1038  \\
$3$  &   0.0354  &   0.0403  &  & \\
$4$  &    0.0193 &  &  & \\
\hline
\end{tabular}
\caption{Test 1. Relative errors $e_{\eta}$ computed for  reconstruction of the function $\eta(t)= 0.7 e^{-t} + 0.05, t \in [0, 300]$ for $T_1= 25, 50, 100$  on different locally adaptively refined time-meshes.}
\label{tabtest1}
\end{table}




\begin{figure}[h!]
\begin{center}
  \begin{tabular}{ccc}
    \hline
    \multicolumn{3}{c}{  $nr.ref.=0$ }\\
  {\includegraphics[scale=0.3, clip=]{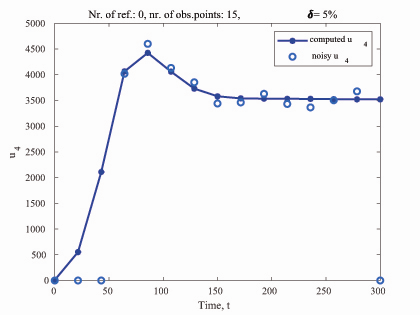}}  &
  {\includegraphics[scale=0.3, clip=]{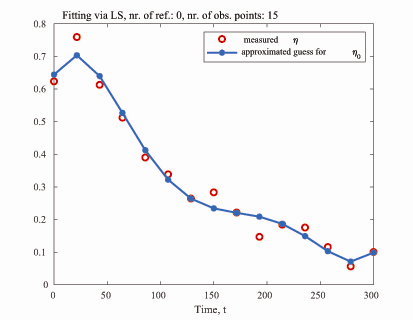}} &
  {\includegraphics[scale=0.3, clip=]{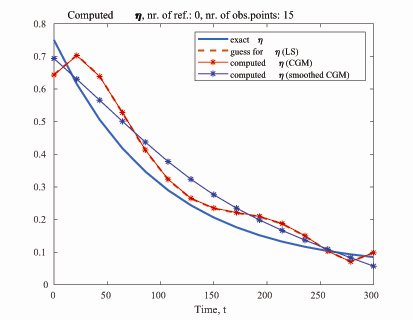}} \\
  ${u_4}_\tau$  & LS fitting to $\eta_\tau(t)$ & $\eta_\tau(t)$  \\
   \hline
  \multicolumn{3}{c}{  $nr.ref.=1$ }\\
  {\includegraphics[scale=0.3, clip=]{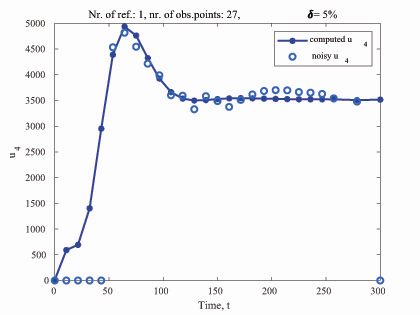}}  &
  {\includegraphics[scale=0.3, clip=]{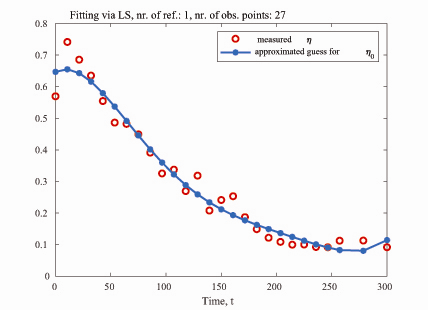}} &
  {\includegraphics[scale=0.3, clip=]{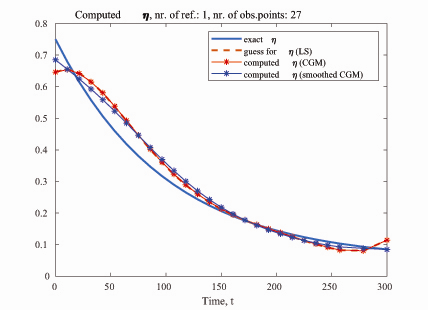}} \\
 ${u_4}_\tau$  & LS fitting to $\eta_\tau(t)$ & $\eta_\tau(t)$  \\
 \hline
   \multicolumn{3}{c}{  $nr.ref.=2$ }\\
  {\includegraphics[scale=0.3, clip=]{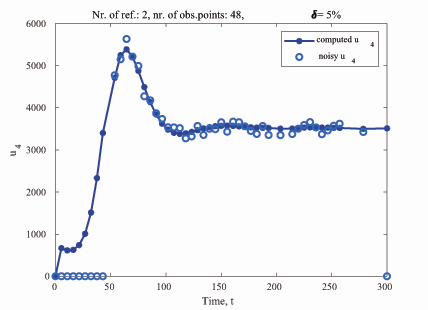}}  &
  {\includegraphics[scale=0.3, clip=]{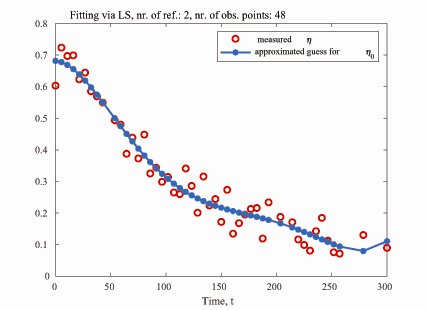}} &
  {\includegraphics[scale=0.3, clip=]{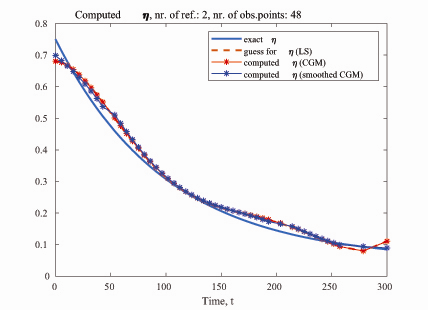}} \\
 ${u_4}_\tau$  & LS fitting to $\eta_\tau(t)$ & $\eta_\tau(t)$  \\
 \hline
    \multicolumn{3}{c}{  $nr.ref.=3$ }\\
  {\includegraphics[scale=0.3, clip=]{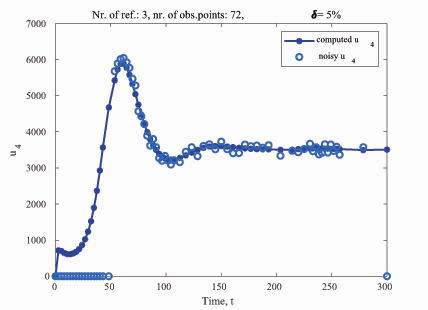}}  &
  {\includegraphics[scale=0.3, clip=]{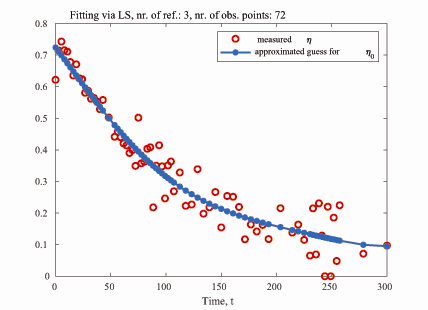}} &
  {\includegraphics[scale=0.3, clip=]{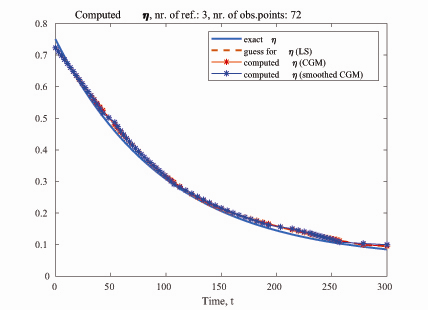}} \\
 ${u_4}_\tau$  & LS fitting to $\eta_\tau(t)$ &  $\eta_\tau(t)$  \\
   \hline
\end{tabular}
\end{center}
\caption{
  \emph{  Test 1. Left figures: simulated ${u_4}_\tau$ vs.  noisy ${u_4}_\tau$ on different adaptively refined time meshes. Here, noisy observed  data are presented by circles. Middle figures: least squares fitting to noisy data for $\eta_\tau$. Right figures: results of ACGA on adaptively refined meshes.  Computations are done for  noise level $\sigma=5\%$ in $u_4$  and for $T_1 = 50$.}}
\label{fig:FIG1}
\end{figure}


\begin{figure}[h!]
\begin{center}
  \begin{tabular}{ccc}
    \hline
    \multicolumn{3}{c}{  $nr.ref.=0$ }\\
  {\includegraphics[scale=0.3, clip=]{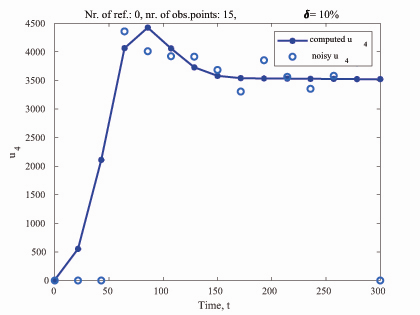}}  &
  {\includegraphics[scale=0.3, clip=]{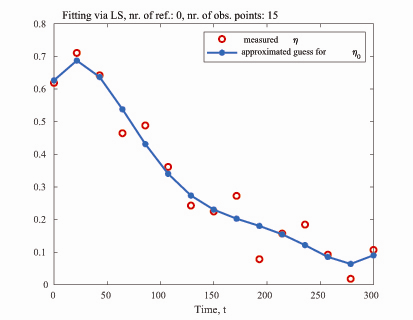}} &
  {\includegraphics[scale=0.3, clip=]{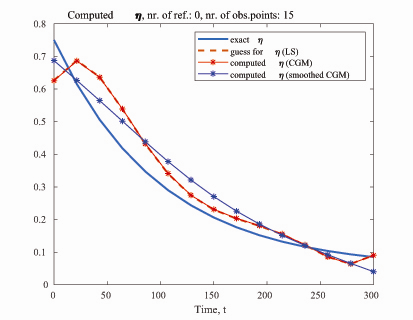}} \\
  ${u_4}_\tau$  & LS fitting to $\eta_\tau(t)$ & $\eta_\tau(t)$  \\
   \hline
  \multicolumn{3}{c}{  $nr.ref.=1$ }\\
  {\includegraphics[scale=0.3, clip=]{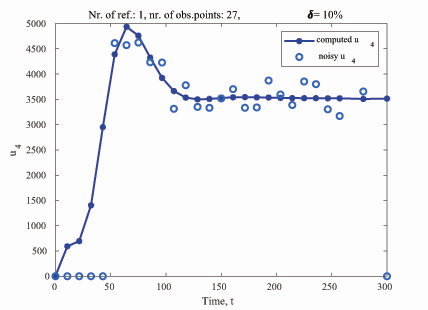}}  &
  {\includegraphics[scale=0.3, clip=]{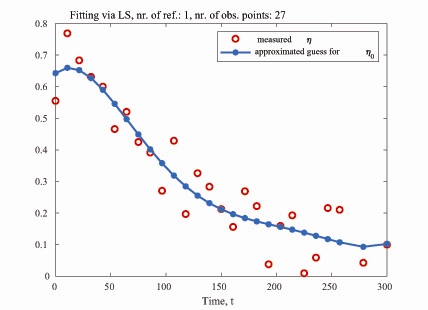}} &
  {\includegraphics[scale=0.3, clip=]{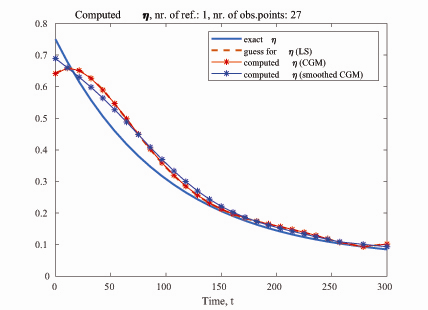}} \\
 ${u_4}_\tau$  & LS fitting to $\eta_\tau(t)$ & $\eta_\tau(t)$  \\
 \hline
   \multicolumn{3}{c}{  $nr.ref.=2$ }\\
  {\includegraphics[scale=0.3, clip=]{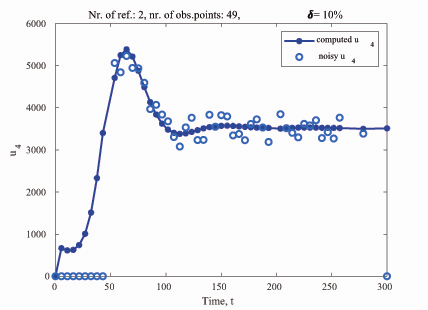}}  &
  {\includegraphics[scale=0.3, clip=]{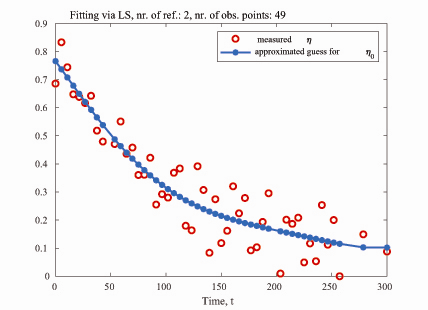}} &
  {\includegraphics[scale=0.3, clip=]{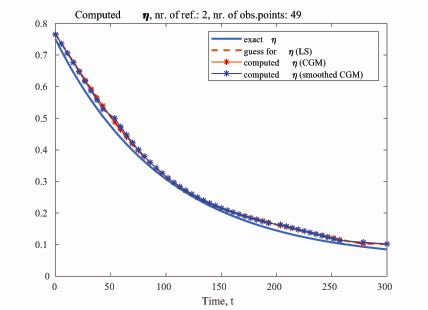}} \\
 ${u_4}_\tau$  & LS fitting to $\eta_\tau(t)$ & $\eta_\tau(t)$  \\
 \hline
    \multicolumn{3}{c}{  $nr.ref.=3$ }\\
  {\includegraphics[scale=0.3, clip=]{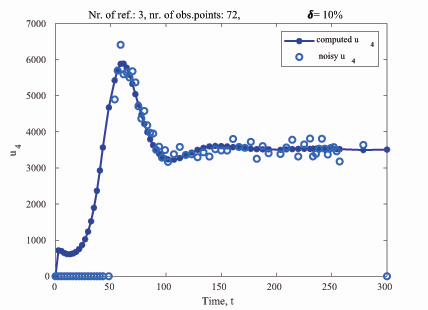}}  &
  {\includegraphics[scale=0.3, clip=]{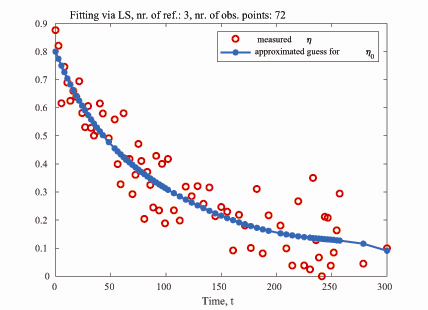}} &
  {\includegraphics[scale=0.3, clip=]{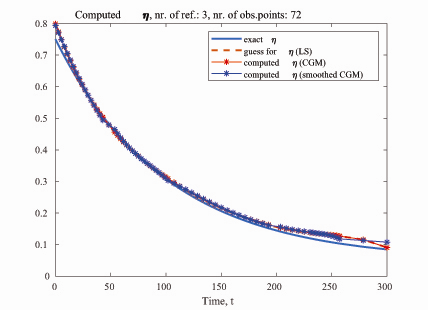}} \\
 ${u_4}_\tau$  & LS fitting to $\eta_\tau(t)$ &  $\eta_\tau(t)$  \\
   \hline
\end{tabular}
\end{center}
\caption{
  \emph{  Test 1. Left figures: simulated ${u_4}_\tau$ vs.  noisy ${u_4}_\tau$ on different adaptively refined time meshes. Here, noisy observed  data are presented by circles. Middle figures: least squares fitting to noisy data for $\eta_\tau$. Right figures: results of ACGA on adaptively refined meshes.  Computations are done for  noise level $\sigma=10\%$ in $u_4$  and for $T_1 = 50$.}}
\label{fig:FIG2}
\end{figure}


\begin{figure}[h!]
\begin{center}
  \begin{tabular}{ccc}
    \hline
    \multicolumn{3}{c}{  $nr.ref.=0$ }\\
  {\includegraphics[scale=0.3, clip=]{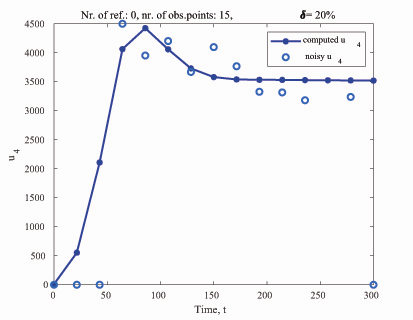}}  &
  {\includegraphics[scale=0.3, clip=]{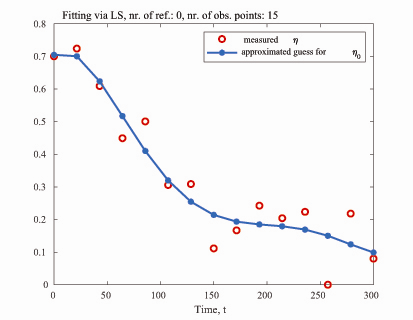}} &
  {\includegraphics[scale=0.3, clip=]{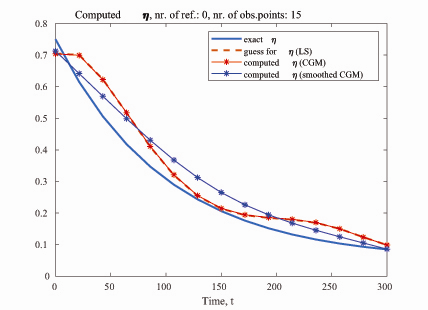}} \\
  ${u_4}_\tau$  & LS fitting to $\eta_\tau(t)$ & $\eta_\tau(t)$  \\
   \hline
  \multicolumn{3}{c}{  $nr.ref.=1$ }\\
  {\includegraphics[scale=0.3, clip=]{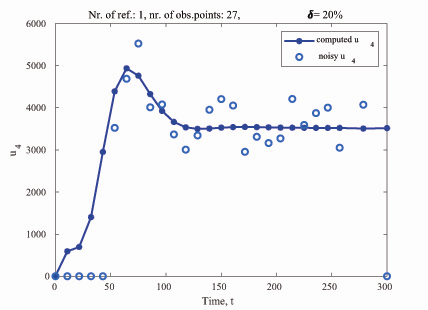}}  &
  {\includegraphics[scale=0.3, clip=]{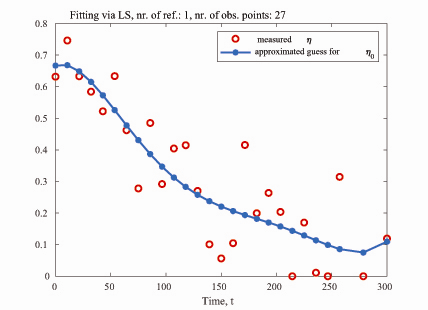}} &
  {\includegraphics[scale=0.3, clip=]{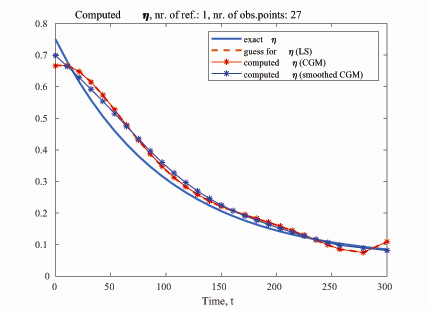}} \\
 ${u_4}_\tau$  & LS fitting to $\eta_\tau(t)$ & $\eta_\tau(t)$  \\
 \hline
   \multicolumn{3}{c}{  $nr.ref.=2$ }\\
  {\includegraphics[scale=0.3, clip=]{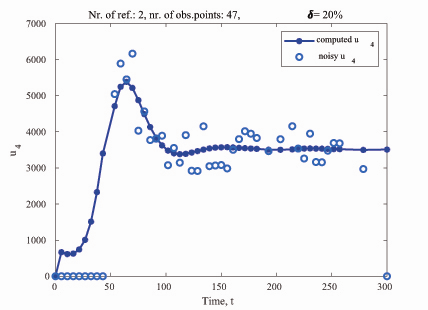}}  &
  {\includegraphics[scale=0.3, clip=]{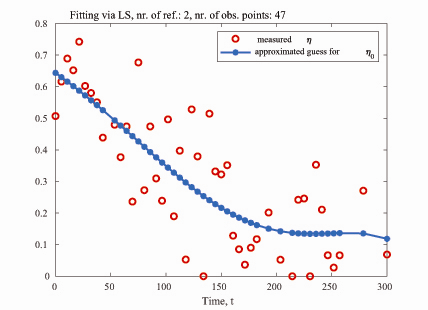}} &
  {\includegraphics[scale=0.3, clip=]{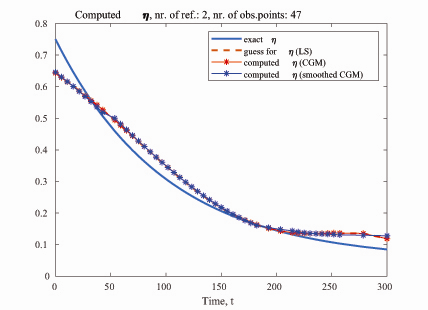}} \\
 ${u_4}_\tau$  & LS fitting to $\eta_\tau(t)$ & $\eta_\tau(t)$  \\
 \hline
    \multicolumn{3}{c}{  $nr.ref.=3$ }\\
  {\includegraphics[scale=0.3, clip=]{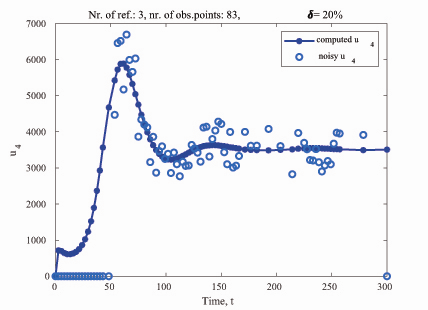}}  &
  {\includegraphics[scale=0.3, clip=]{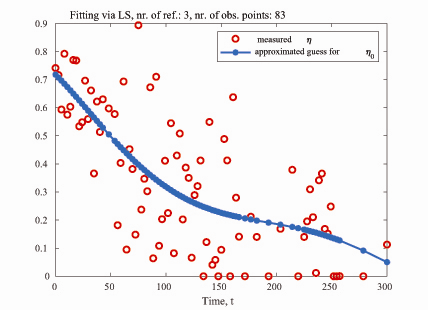}} &
  {\includegraphics[scale=0.3, clip=]{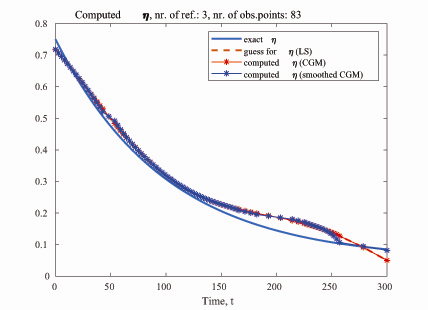}} \\
 ${u_4}_\tau$  & LS fitting to $\eta_\tau(t)$ &  $\eta_\tau(t)$  \\
   \hline
\end{tabular}
\end{center}
\caption{
  \emph{  Test 1. Left figures: simulated ${u_4}_\tau$ vs.  noisy ${u_4}_\tau$ on different adaptively refined time meshes. Here, noisy observed  data are presented by circles. Middle figures: least squares fitting to noisy data for $\eta_\tau$. Right figures: results of ACGA on adaptively refined meshes.  Computations are done for  noise level $\sigma=20\%$ in $u_4$  and for $T_1 = 50$.}}
\label{fig:FIG3}
\end{figure}


\begin{figure}[h!]
\begin{center}
  \begin{tabular}{ccc}
    \hline
    \multicolumn{3}{c}{  $nr.ref.=0$ }\\
   {\includegraphics[scale=0.3, clip=]{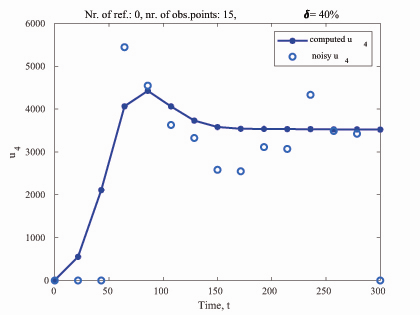}}  &
  {\includegraphics[scale=0.3, clip=]{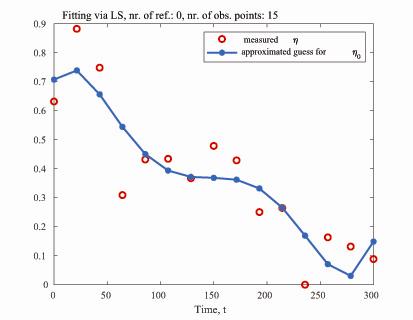}} &
  {\includegraphics[scale=0.3, clip=]{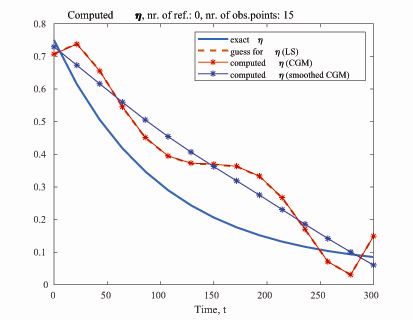}} \\
  ${u_4}_\tau$  & LS fitting to $\eta_\tau(t)$ & $\eta_\tau(t)$  \\
   \hline
  \multicolumn{3}{c}{  $nr.ref.=1$ }\\
  {\includegraphics[scale=0.3, clip=]{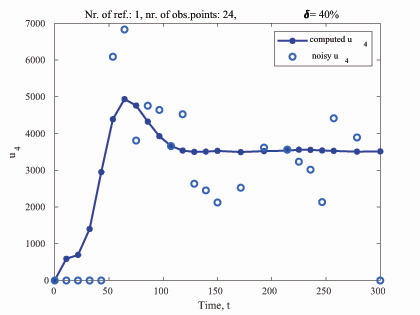}}  &
  {\includegraphics[scale=0.3, clip=]{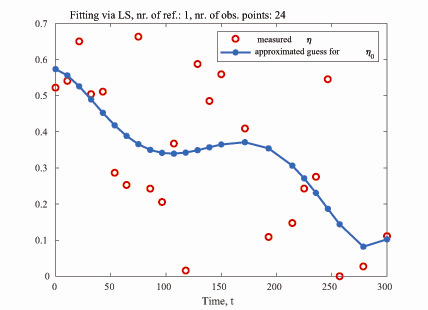}} &
  {\includegraphics[scale=0.3, clip=]{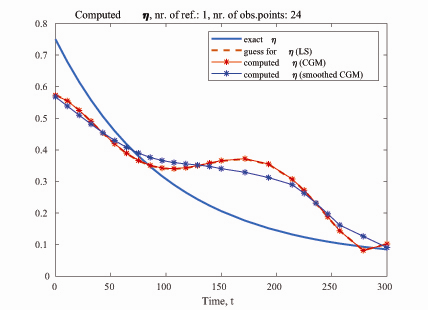}} \\
 ${u_4}_\tau$  & LS fitting to $\eta_\tau(t)$ & $\eta_\tau(t)$  \\
 \hline
   \multicolumn{3}{c}{  $nr.ref.=2$ }\\
  {\includegraphics[scale=0.3, clip=]{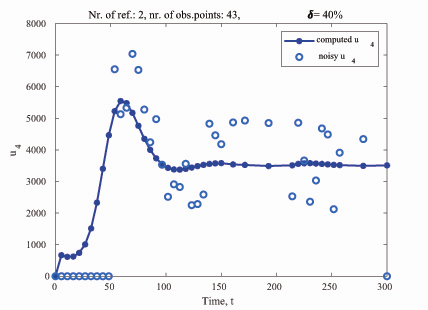}}  &
  {\includegraphics[scale=0.3, clip=]{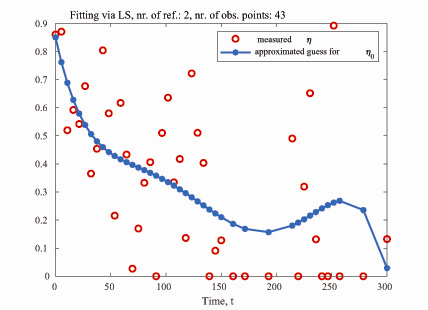}} &
  {\includegraphics[scale=0.3, clip=]{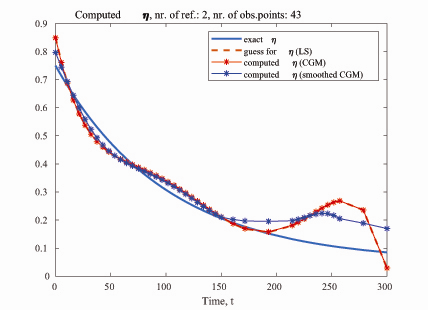}} \\
 ${u_4}_\tau$  & LS fitting to $\eta_\tau(t)$ & $\eta_\tau(t)$  \\
 \hline
    \multicolumn{3}{c}{  $nr.ref.=3$ }\\
  {\includegraphics[scale=0.3, clip=]{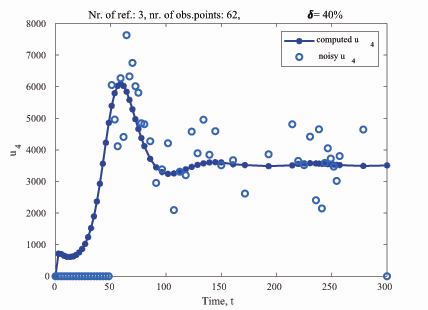}}  &
  {\includegraphics[scale=0.3, clip=]{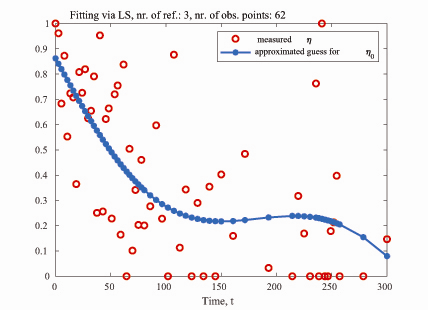}} &
  {\includegraphics[scale=0.3, clip=]{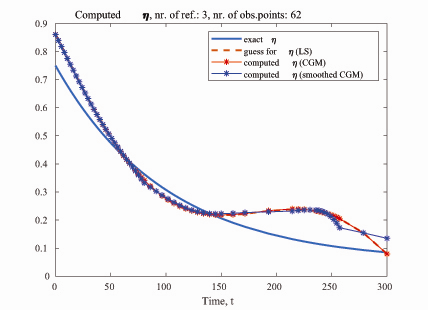}} \\
 ${u_4}_\tau$  & LS fitting to $\eta_\tau(t)$ &  $\eta_\tau(t)$  \\
   \hline
\end{tabular}
\end{center}
\caption{
  \emph{  Test 1. Left figures: simulated ${u_4}_\tau$ vs.  noisy ${u_4}_\tau$ on different adaptively refined time meshes. Here, noisy observed  data are presented by circles. Middle figures: least squares fitting to noisy data for $\eta_\tau$. Right figures: results of ACGA on adaptively refined meshes.  Computations are done for  noise level $\sigma=40\%$ in $u_4$  and for $T_1 = 50$.}}
\label{fig:FIG4}
\end{figure}

\begin{figure}[h!]
\begin{center}
  \begin{tabular}{ccc}
    \hline
    \multicolumn{3}{c}{  $nr.ref.=0$ }\\
  {\includegraphics[scale=0.3, clip=]{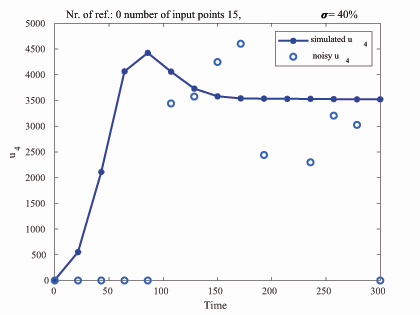}}  &
  {\includegraphics[scale=0.3, clip=]{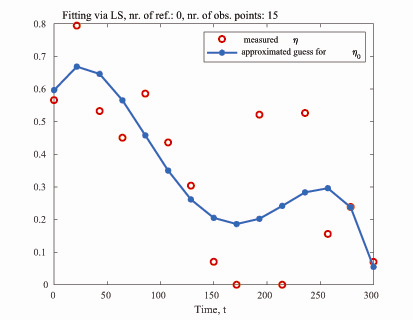}} &
  {\includegraphics[scale=0.3, clip=]{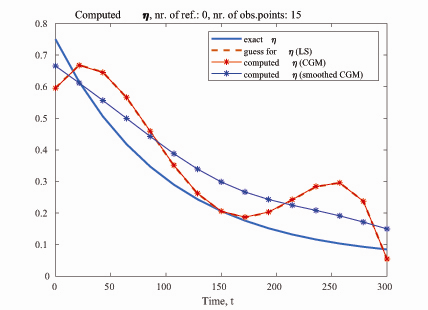}} \\
  ${u_4}_\tau$  & LS fitting to $\eta_\tau(t)$ & $\eta_\tau(t)$  \\
   \hline
  \multicolumn{3}{c}{  $nr.ref.=1$ }\\
  {\includegraphics[scale=0.3, clip=]{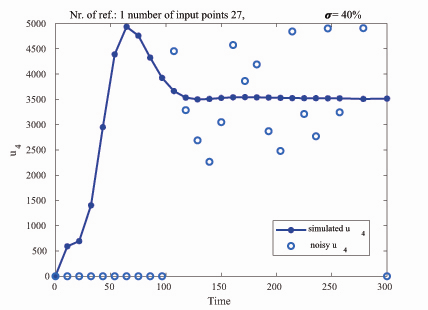}}  &
  {\includegraphics[scale=0.3, clip=]{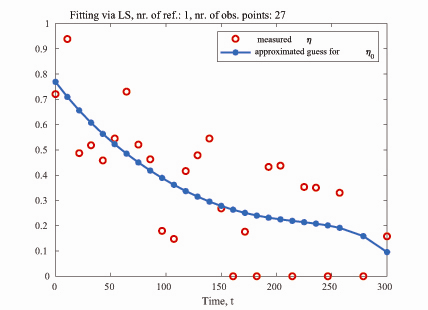}} &
  {\includegraphics[scale=0.3, clip=]{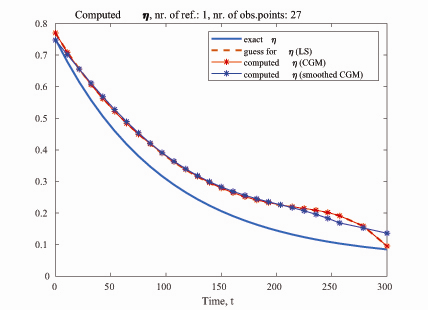}} \\
 ${u_4}_\tau$  & LS fitting to $\eta_\tau(t)$ & $\eta_\tau(t)$  \\
 \hline
   \multicolumn{3}{c}{  $nr.ref.=2$ }\\
  {\includegraphics[scale=0.3, clip=]{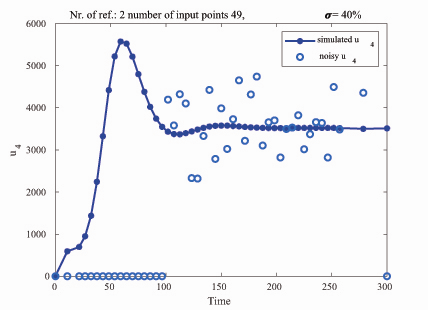}}  &
  {\includegraphics[scale=0.3, clip=]{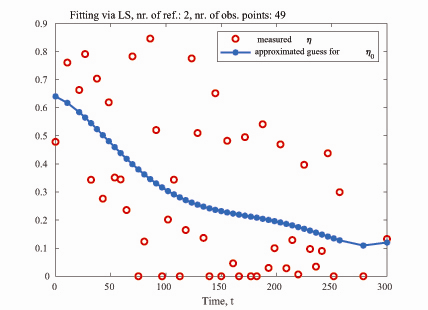}} &
  {\includegraphics[scale=0.3, clip=]{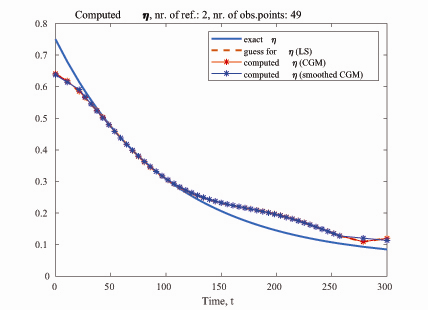}} \\
 ${u_4}_\tau$  & LS fitting to $\eta_\tau(t)$ & $\eta_\tau(t)$  \\
   \hline
\end{tabular}
\end{center}
\caption{
  \emph{  Test 1. Left figures: simulated ${u_4}_\tau$ vs.  noisy ${u_4}_\tau$ on different adaptively refined time meshes. Here, noisy observed  data are presented by circles.  Middle figures: least squares fitting to noisy data for $\eta_\tau$. Right figures: results of ACGA on adaptively refined meshes.  Computations are done for  noise level $\sigma=40\%$ in $u_4$  and for $T_1 = 100$.}}
\label{fig:FIG5}
\end{figure}


In this test we present results of reconstruction of a smooth function
$\eta(t)= 0.7 e^{-t} + 0.05, t \in [0, 300]$ for $T_1= 25, 50, 100$ and number of observation
 points $15$. The initial time partition $\mathcal J_{\tau}$
 is generated with equidistant time step $\tau = 300/14$.  Results of
 reconstruction of the model function $\eta(t)= 0.7 e^{-t} + 0.05$ for
 noise levels $\sigma = 5\%, 10\%, 20\%, 40\%$ in data $u_4(t)$ are
 presented Table 1.  Figures \ref{fig:FIG1}-\ref{fig:FIG4} show results of
 reconstruction of the function $\eta(t)= 0.7 e^{-t} + 0.05$ for noise
 levels $\sigma = 5\%, 10\%, 20\%, 40\%$ in data $u_4(t)$ for $T_1 = 50$, respectively. Figure \ref{fig:FIG5} shows results of reconstruction of this function for  noise level $\sigma = 40\%$ in data $u_4(t)$ and for $T_1 = 100$.

 Table 1 and Figures \ref{fig:FIG1}-\ref{fig:FIG5} confirm that with local time-mesh refinements the reconstruction of the drug efficacy function $\eta_\tau$ is significantly improved compared to the reconstruction of $\eta_\tau$ obtained on initial non-refined time-mesh.


\subsection{Test 2}

\begin{table}[h!]
\center
\begin{tabular}{| l|l|l|l|l| }
  \multicolumn{5}{c}{  $T_1 = 25$ }\\
  \hline
$\sigma$ & 5 \%  &   10\%   &  20\%   & 40\%       \\
nr.of ref. &   &    &   &        \\
\hline
$0$ &    0.0718   &  0.0802   &  0.0834  &      0.0617               \\
$1$  &  0.0592  &    0.0315   &  0.0290   &  0.0493  \\
$2$  &   0.0403  &  0.0091  &  & 0.0301 \\
$3$  &   0.0272  &  0.0050  &  &  0.0240 \\
$4$  &   0.0191  &  0.0064 &  & \\
$5$  & 0.0170  &  &  & \\
$6$  &   0.0117  &  &  & \\
\hline
\multicolumn{5}{c}{  $T_1 = 50$ }\\
  \hline
$\sigma$ & 5 \%  &   10\%   &  20\%   & 40\%       \\
nr.of ref. &   &    &   &        \\
\hline
$0$ &    0.0725     &   0.0758    &  0.0720  &  0.1026                  \\
$1$  &   0.0656 &  0.0572 &   0.0694   &  0.0730 \\
$2$  &    0.0459   &  0.0414 &  0.0505 &  0.0571  \\
$3$  &   0.0273  &  0.0239   &  0.0179 &  0.0236\\
$4$  &  0.0111  & 0.0183  &  & \\
$5$  &   0.0066    &  0.0099  &  & \\
\hline
\multicolumn{5}{c}{  $T_1 = 100$ }\\
  \hline
$\sigma$ & 5 \%  &   10\%   &  20\%   & 40\%       \\
nr.of ref. &   &    &   &        \\
\hline
$0$ &   0.0801  &  0.0676   &   0.0535 &     0.0852               \\
$1$  &   0.0568   & 0.0547  &  &  0.0487  \\
$2$  &   0.0351  &   0.0481   &  & 0.0208  \\
$3$  & 0.0265   &  0.0265 &  & \\
$4$  & 0.0212   & 0.0130  &  & \\
$5$  &  0.0095   & 0.0090  &  & \\
$6$  &   0.0084     &  &  & \\
\hline
\end{tabular}
\caption{
  Test 2. Relative errors  $e_{\eta}$ computed for  reconstruction of the function $\eta(t)= 0.7, t \in [0, 300]$ for $T_1= 25, 50, 100$  on different locally adaptively refined time-meshes.}
\label{tabtest2}
\end{table}



\begin{figure}[h!]
\begin{center}
  \begin{tabular}{ccc}
    \hline
    \multicolumn{3}{c}{  $nr.ref.=0$ }\\
  {\includegraphics[scale=0.3, clip=]{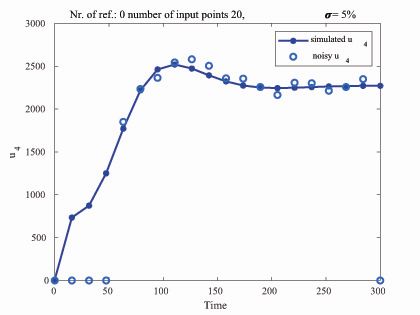}} &
  {\includegraphics[scale=0.3, clip=]{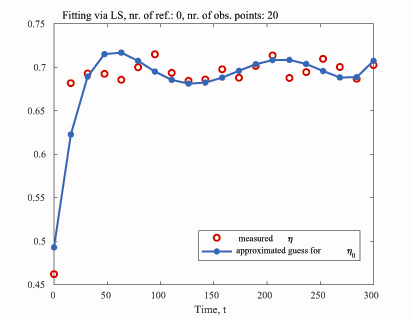}} &
  {\includegraphics[scale=0.3, clip=]{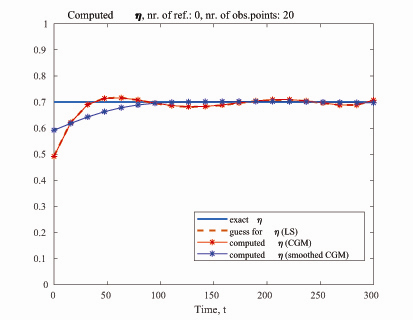}} \\
  ${u_4}_\tau$  & LS fitting to $\eta_\tau(t)$ & $\eta_\tau(t)$  \\
   \hline
  \multicolumn{3}{c}{  $nr.ref.=1$ }\\
  {\includegraphics[scale=0.3, clip=]{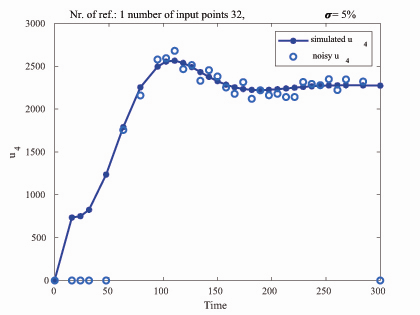}}  &
  {\includegraphics[scale=0.3, clip=]{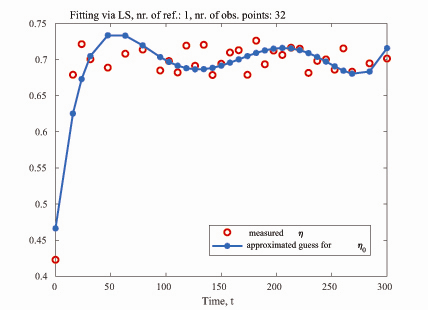}} &
  {\includegraphics[scale=0.3, clip=]{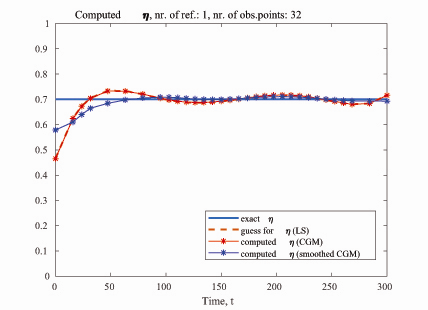}} \\
 ${u_4}_\tau$  & LS fitting to $\eta_\tau(t)$ & $\eta_\tau(t)$  \\
 \hline
   \multicolumn{3}{c}{  $nr.ref.=2$ }\\
  {\includegraphics[scale=0.3, clip=]{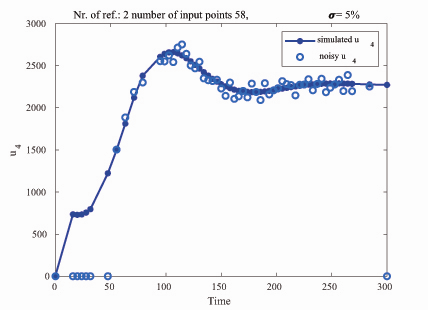}}  &
  {\includegraphics[scale=0.3, clip=]{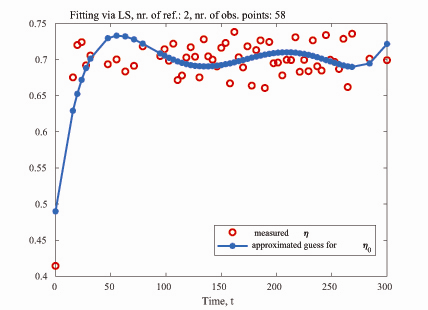}} &
  {\includegraphics[scale=0.3, clip=]{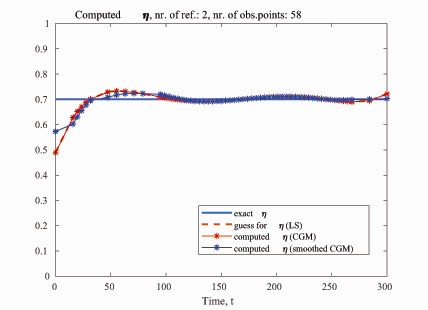}} \\
 ${u_4}_\tau$  & LS fitting to $\eta_\tau(t)$ & $\eta_\tau(t)$  \\
 \hline
   \multicolumn{3}{c}{  $nr.ref.=3$ }\\
  {\includegraphics[scale=0.3, clip=]{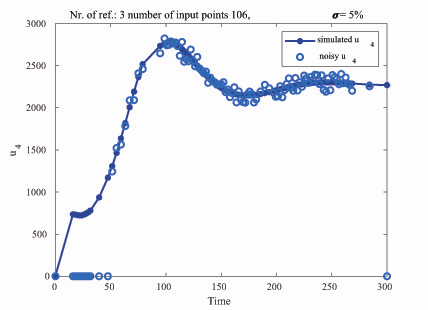}}  &
  {\includegraphics[scale=0.3, clip=]{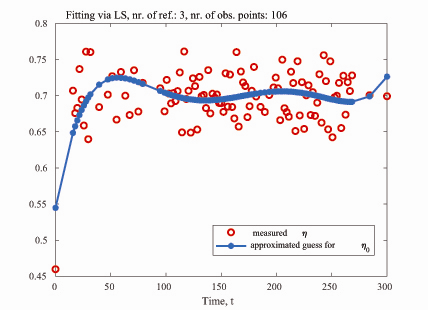}} &
  {\includegraphics[scale=0.3, clip=]{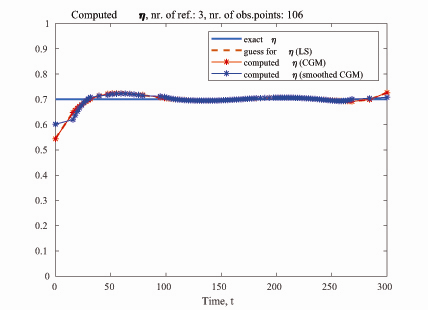}} \\
 ${u_4}_\tau$  & LS fitting to $\eta_\tau(t)$ & $\eta_\tau(t)$  \\
   \hline
\end{tabular}
\end{center}
\caption{
  \emph{  Test 2.  Left figures: simulated ${u_4}_\tau$ vs.  noisy ${u_4}_\tau$ on different adaptively refined time meshes. Here, noisy observed  data are presented by circles.  Middle figures: least squares fitting to noisy data for $\eta_\tau$. Right figures: results of ACGA on adaptively refined meshes.  Computations are done for  noise level $\sigma=5\%$ in $u_4$  and for $T_1 = 50$.}}
\label{fig:FIG6}
\end{figure}


\begin{figure}[h!]
\begin{center}
  \begin{tabular}{ccc}
    \hline
    \multicolumn{3}{c}{  $nr.ref.=0$ }\\
  {\includegraphics[scale=0.3, clip=]{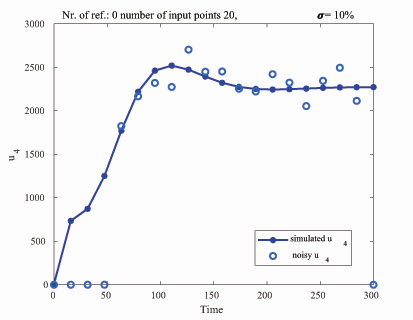}}  &
  {\includegraphics[scale=0.3, clip=]{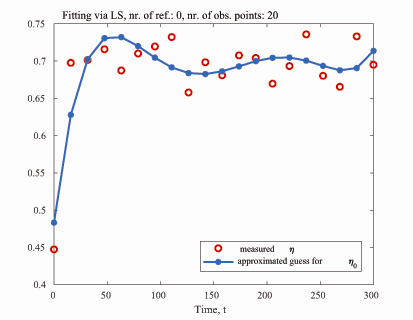}} &
  {\includegraphics[scale=0.3, clip=]{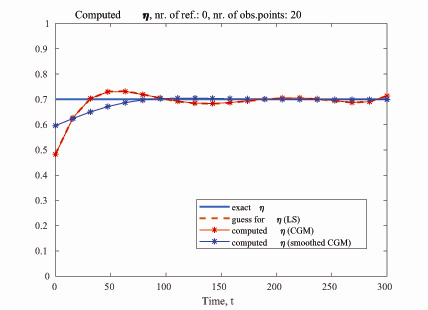}} \\
  ${u_4}_\tau$  & LS fitting to $\eta_\tau(t)$ & $\eta_\tau(t)$  \\
   \hline
  \multicolumn{3}{c}{  $nr.ref.=1$ }\\
  {\includegraphics[scale=0.3, clip=]{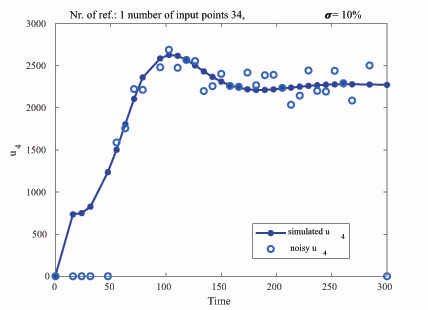}}  &
  {\includegraphics[scale=0.3, clip=]{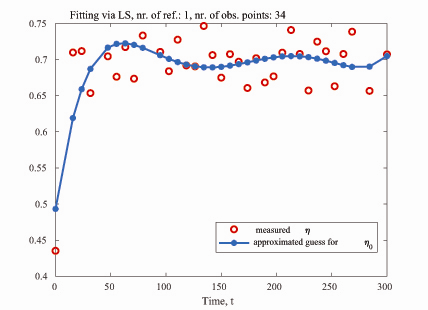}} &
  {\includegraphics[scale=0.3, clip=]{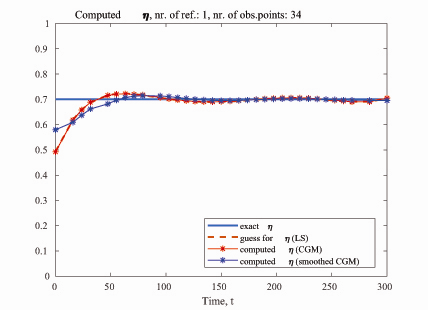}} \\
 ${u_4}_\tau$  & LS fitting to $\eta_\tau(t)$ & $\eta_\tau(t)$  \\
 \hline
   \multicolumn{3}{c}{  $nr.ref.=2$ }\\
  {\includegraphics[scale=0.3, clip=]{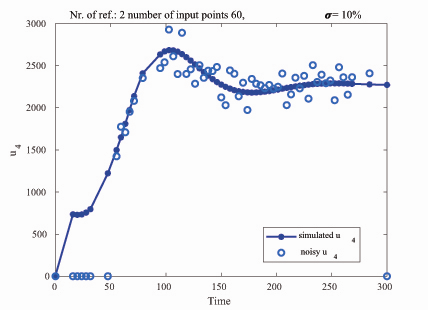}}  &
  {\includegraphics[scale=0.3, clip=]{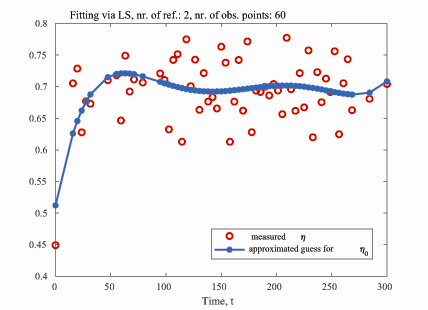}} &
  {\includegraphics[scale=0.3, clip=]{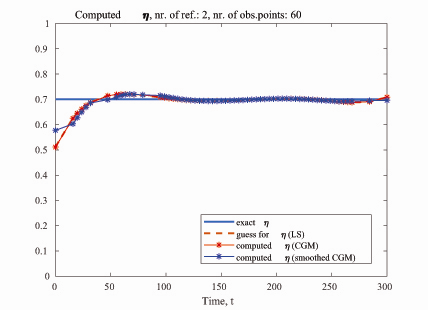}} \\
 ${u_4}_\tau$  & LS fitting to $\eta_\tau(t)$ & $\eta_\tau(t)$  \\
 \hline
   \multicolumn{3}{c}{  $nr.ref.=3$ }\\
  {\includegraphics[scale=0.3, clip=]{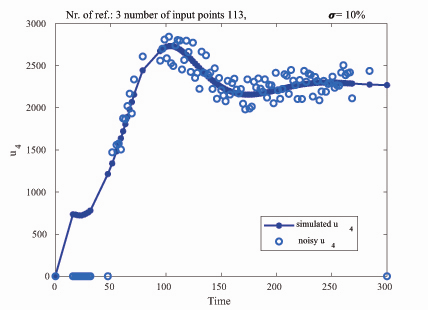}}  &
  {\includegraphics[scale=0.3, clip=]{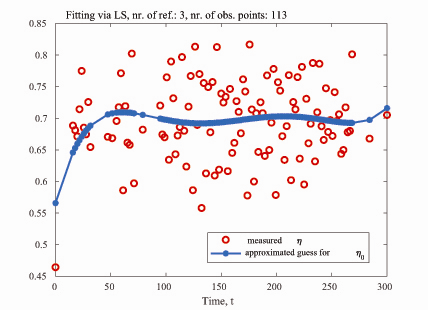}} &
  {\includegraphics[scale=0.3, clip=]{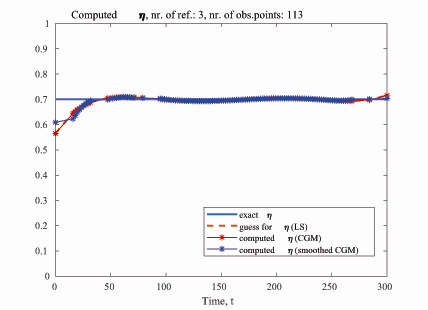}} \\
 ${u_4}_\tau$  & LS fitting to $\eta_\tau(t)$ & $\eta_\tau(t)$  \\
   \hline
\end{tabular}
\end{center}
\caption{
  \emph{  Test 2.  Left figures: simulated ${u_4}_\tau$ vs.  noisy ${u_4}_\tau$ on different adaptively refined time meshes. Here, noisy observed  data are presented by circles.  Middle figures: least squares fitting to noisy data for $\eta_\tau$. Right figures: results of ACGA on adaptively refined meshes.  Computations are done for  noise level $\sigma=10\%$ in $u_4$  and for $T_1 = 50$.}}
\label{fig:FIG7}
\end{figure}


\begin{figure}[h!]
\begin{center}
  \begin{tabular}{ccc}
    \hline
    \multicolumn{3}{c}{  $nr.ref.=0$ }\\
  {\includegraphics[scale=0.3, clip=]{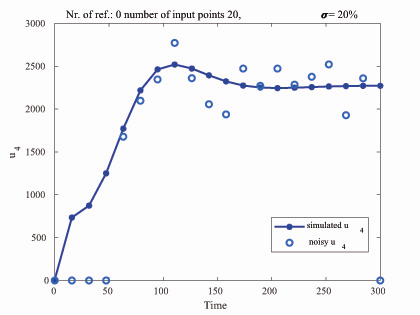}}  &
  {\includegraphics[scale=0.3, clip=]{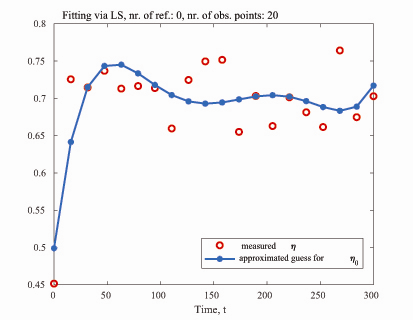}} &
  {\includegraphics[scale=0.3, clip=]{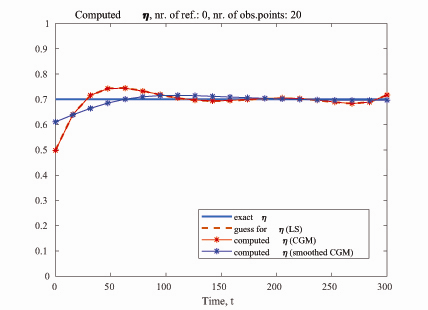}} \\
  ${u_4}_\tau$  & LS fitting to $\eta_\tau(t)$ & $\eta_\tau(t)$  \\
   \hline
  \multicolumn{3}{c}{  $nr.ref.=1$ }\\
  {\includegraphics[scale=0.3, clip=]{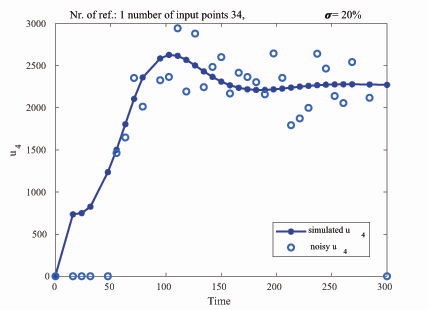}}  &
  {\includegraphics[scale=0.3, clip=]{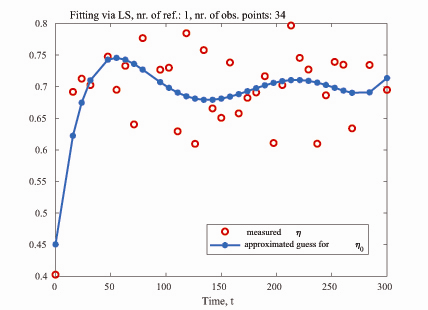}} &
  {\includegraphics[scale=0.3, clip=]{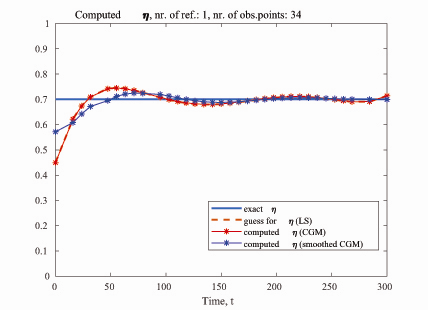}} \\
 ${u_4}_\tau$  & LS fitting to $\eta_\tau(t)$ & $\eta_\tau(t)$  \\
 \hline
   \multicolumn{3}{c}{  $nr.ref.=2$ }\\
  {\includegraphics[scale=0.3, clip=]{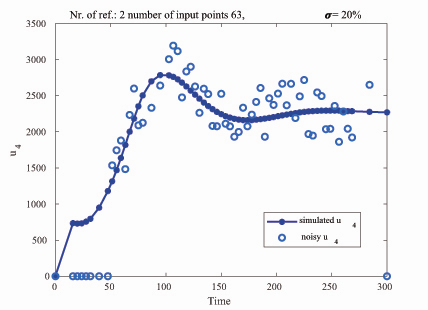}}  &
  {\includegraphics[scale=0.3, clip=]{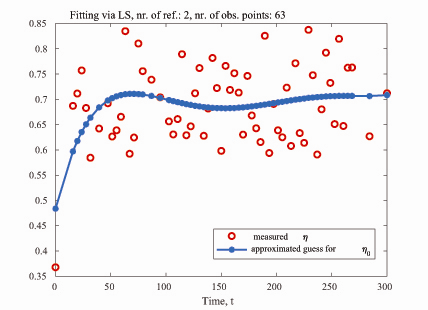}} &
  {\includegraphics[scale=0.3, clip=]{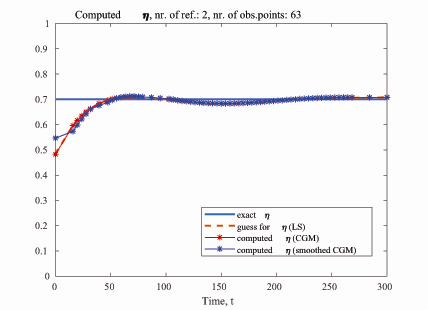}} \\
 ${u_4}_\tau$  & LS fitting to $\eta_\tau(t)$ & $\eta_\tau(t)$  \\
 \hline
   \multicolumn{3}{c}{  $nr.ref.=3$ }\\
  {\includegraphics[scale=0.3, clip=]{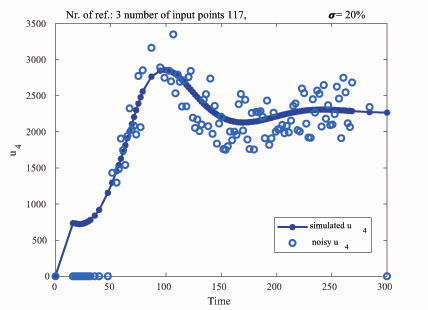}}  &
  {\includegraphics[scale=0.3, clip=]{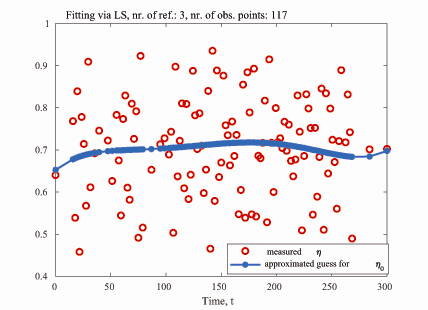}} &
  {\includegraphics[scale=0.3, clip=]{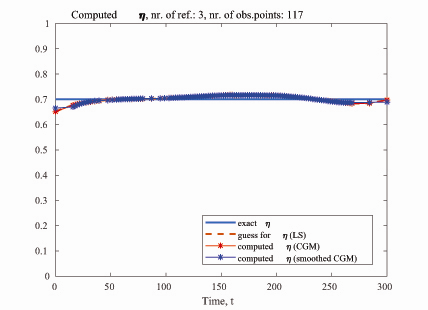}} \\
 ${u_4}_\tau$  & LS fitting to $\eta_\tau(t)$ & $\eta_\tau(t)$  \\
   \hline
\end{tabular}
\end{center}
\caption{
  \emph{  Test 2.  Left figures: simulated ${u_4}_\tau$ vs.  noisy ${u_4}_\tau$ on different adaptively refined time meshes. Here, noisy observed  data are presented by circles.  Middle figures: least squares fitting to noisy data for $\eta_\tau$. Right figures: results of ACGA on adaptively refined meshes.  Computations are done for  noise level $\sigma=20\%$ in $u_4$  and for $T_1 = 50$.}}
\label{fig:FIG8}
\end{figure}


\begin{figure}[h!]
\begin{center}
  \begin{tabular}{ccc}
    \hline
    \multicolumn{3}{c}{  $nr.ref.=0$ }\\
  {\includegraphics[scale=0.3, clip=]{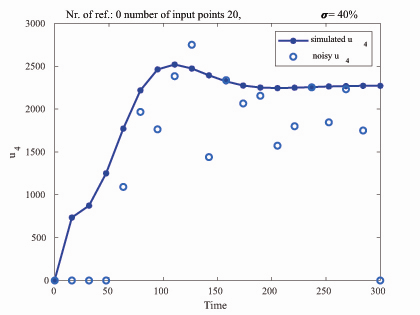}}  &
  {\includegraphics[scale=0.3, clip=]{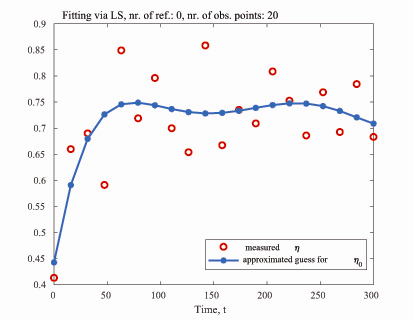}} &
  {\includegraphics[scale=0.3, clip=]{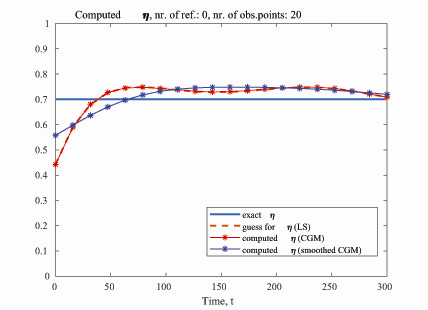}} \\
  ${u_4}_\tau$  & LS fitting to $\eta_\tau(t)$ & $\eta_\tau(t)$  \\
   \hline
  \multicolumn{3}{c}{  $nr.ref.=1$ }\\
  {\includegraphics[scale=0.3, clip=]{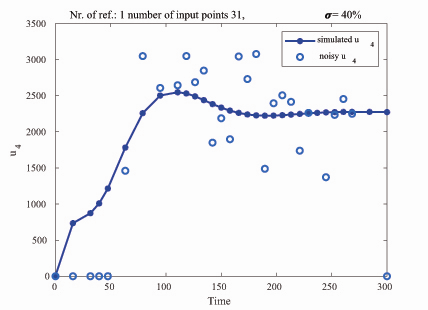}}  &
  {\includegraphics[scale=0.3, clip=]{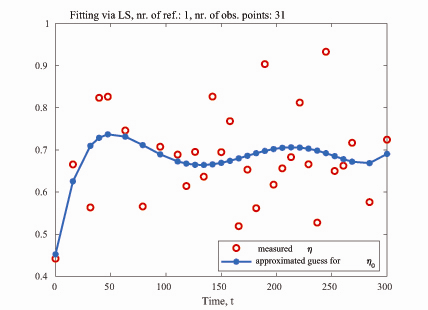}} &
  {\includegraphics[scale=0.3, clip=]{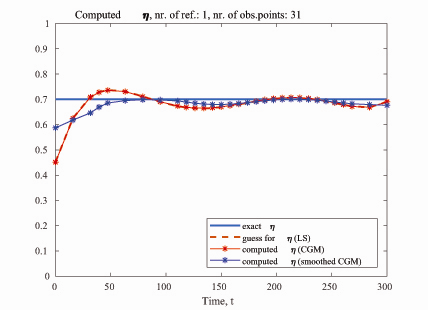}} \\
 ${u_4}_\tau$  & LS fitting to $\eta_\tau(t)$ & $\eta_\tau(t)$  \\
 \hline
   \multicolumn{3}{c}{  $nr.ref.=2$ }\\
  {\includegraphics[scale=0.3, clip=]{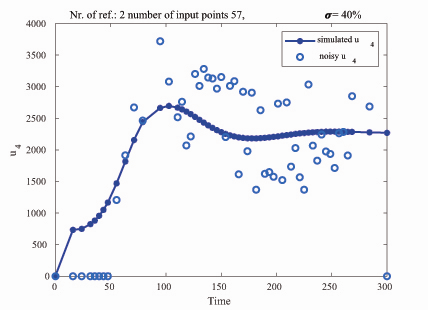}}  &
  {\includegraphics[scale=0.3, clip=]{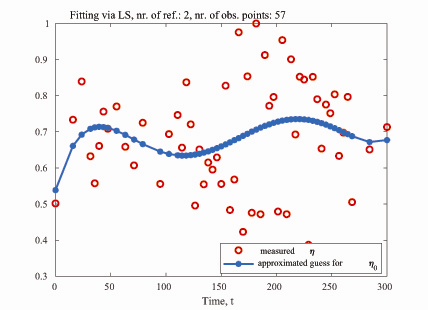}} &
  {\includegraphics[scale=0.3, clip=]{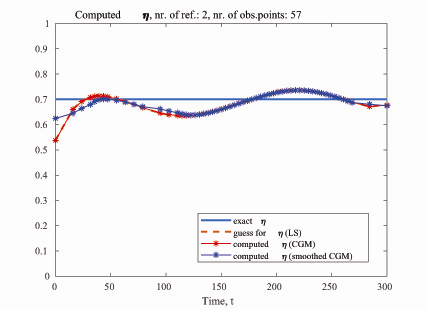}} \\
 ${u_4}_\tau$  & LS fitting to $\eta_\tau(t)$ & $\eta_\tau(t)$  \\
 \hline
   \multicolumn{3}{c}{  $nr.ref.=3$ }\\
  {\includegraphics[scale=0.3, clip=]{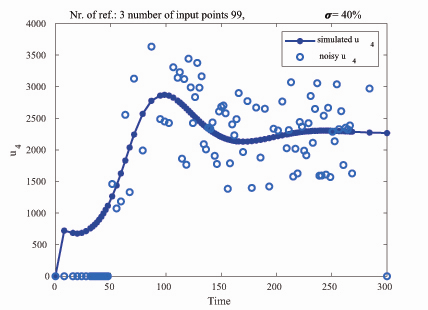}}  &
  {\includegraphics[scale=0.3, clip=]{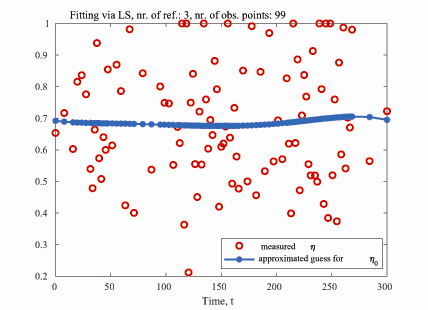}} &
  {\includegraphics[scale=0.3, clip=]{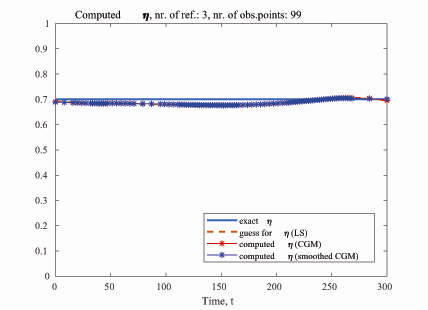}} \\
 ${u_4}_\tau$  & LS fitting to $\eta_\tau(t)$ & $\eta_\tau(t)$  \\
   \hline
\end{tabular}
\end{center}
\caption{
  \emph{  Test 2.  Left figures: simulated ${u_4}_\tau$ vs.  noisy ${u_4}_\tau$ on different adaptively refined time meshes. Here, noisy observed  data are presented by circles.  Middle figures: least squares fitting to noisy data for $\eta_\tau$. Right figures: results of ACGA on adaptively refined meshes.  Computations are done for  noise level $\sigma=40\%$ in $u_4$  and for $T_1 = 50$.}}
\label{fig:FIG9}
\end{figure}


\begin{figure}[h!]
\begin{center}
  \begin{tabular}{ccc}
    \hline
    \multicolumn{3}{c}{  $nr.ref.=0$ }\\
  {\includegraphics[scale=0.3, clip=]{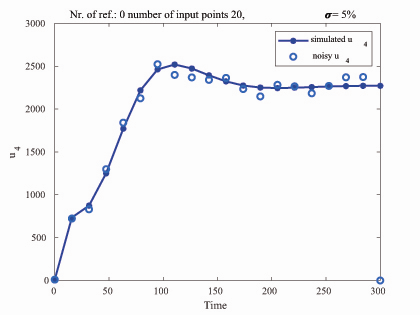}}  &
  {\includegraphics[scale=0.3, clip=]{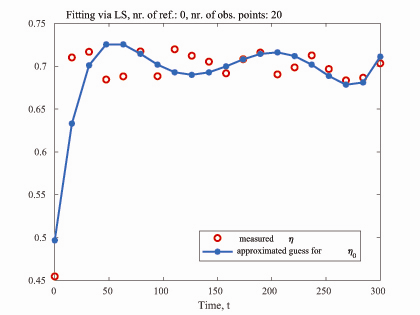}} &
  {\includegraphics[scale=0.3, clip=]{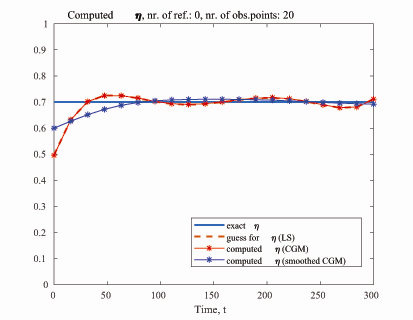}} \\
  ${u_4}_\tau$   & LS fitting to $\eta_\tau(t)$ & $\eta_\tau(t)$  \\
   \hline
  \multicolumn{3}{c}{  $nr.ref.=1$ }\\
  {\includegraphics[scale=0.3, clip=]{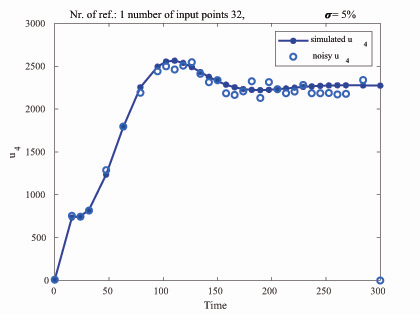}}  &
  {\includegraphics[scale=0.3, clip=]{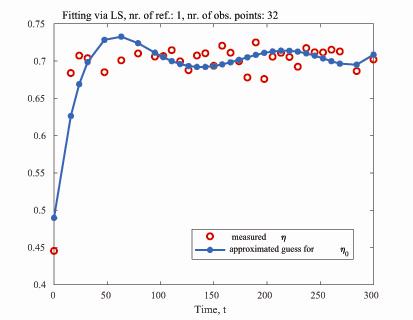}} &
  {\includegraphics[scale=0.3, clip=]{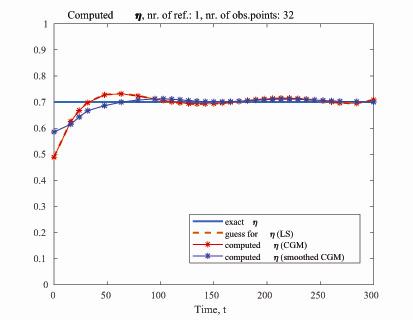}} \\
 ${u_4}_\tau$  & LS fitting to $\eta_\tau(t)$ & $\eta_\tau(t)$  \\
 \hline
   \multicolumn{3}{c}{  $nr.ref.=2$ }\\
  {\includegraphics[scale=0.3, clip=]{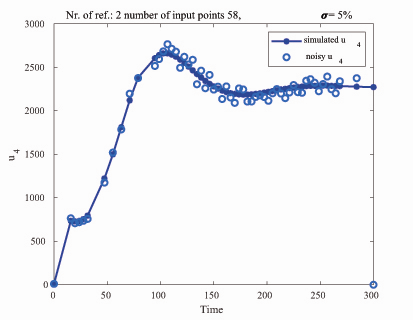}}  &
  {\includegraphics[scale=0.3, clip=]{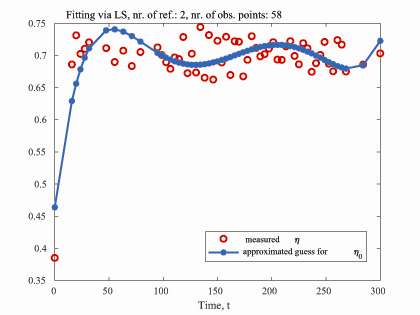}} &
  {\includegraphics[scale=0.3, clip=]{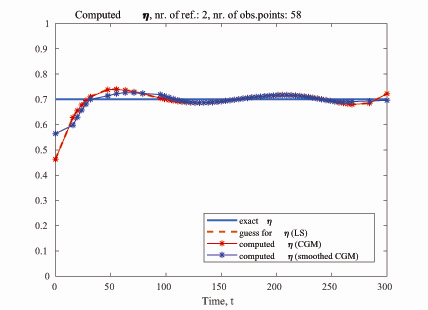}} \\
 ${u_4}_\tau$  & LS fitting to $\eta_\tau(t)$ & $\eta_\tau(t)$  \\
 \hline
\end{tabular}
\end{center}
\caption{
  \emph{  Test 2. Left figures: simulated ${u_4}_\tau$ vs.  noisy ${u_4}_\tau$ on different adaptively refined time meshes. Here, noisy observed  data are presented by circles.   Middle figures: least squares fitting to noisy data for $\eta_\tau$. Right figures: results of ACGA on adaptively refined meshes.  Computations are done for  noise level $\sigma=5\%$ in $u_4$  and for $T_1 = 100$.}}
\label{fig:FIG10}
\end{figure}



\begin{figure}[h!]
\begin{center}
  \begin{tabular}{ccc}
    \hline
    \multicolumn{3}{c}{  $nr.ref.=0$ }\\
  {\includegraphics[scale=0.3, clip=]{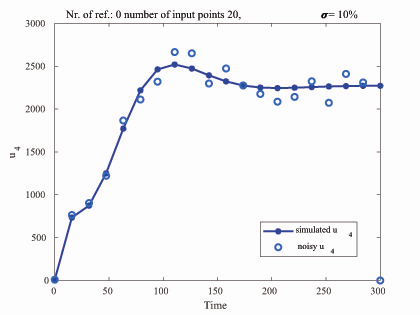}}  &
  {\includegraphics[scale=0.3, clip=]{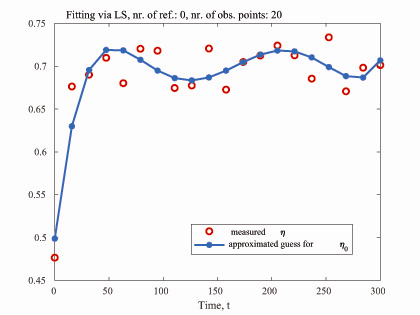}} &
  {\includegraphics[scale=0.3, clip=]{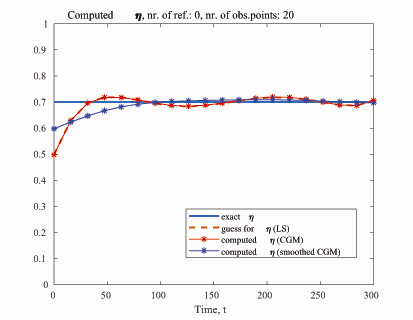}} \\
  ${u_4}_\tau$   & LS fitting to $\eta_\tau(t)$ & $\eta_\tau(t)$  \\
   \hline
  \multicolumn{3}{c}{  $nr.ref.=1$ }\\
  {\includegraphics[scale=0.3, clip=]{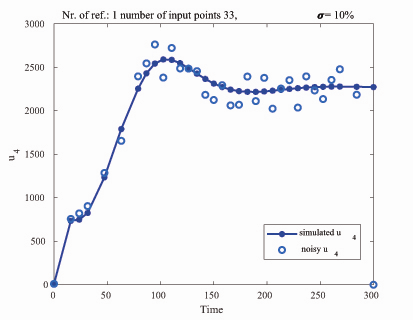}}  &
  {\includegraphics[scale=0.3, clip=]{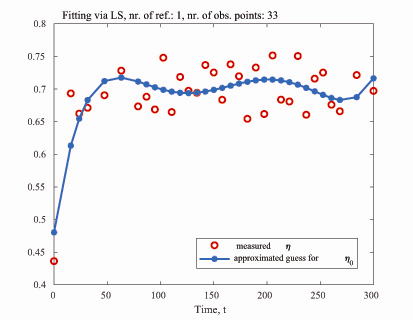}} &
  {\includegraphics[scale=0.3, clip=]{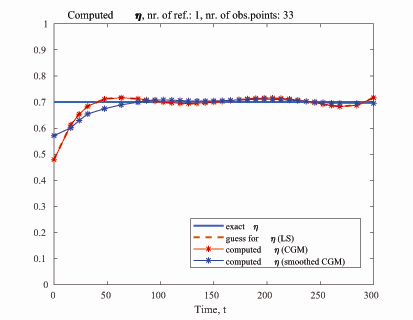}} \\
 ${u_4}_\tau$  & LS fitting to $\eta_\tau(t)$ & $\eta_\tau(t)$  \\
 \hline
   \multicolumn{3}{c}{  $nr.ref.=2$ }\\
  {\includegraphics[scale=0.3, clip=]{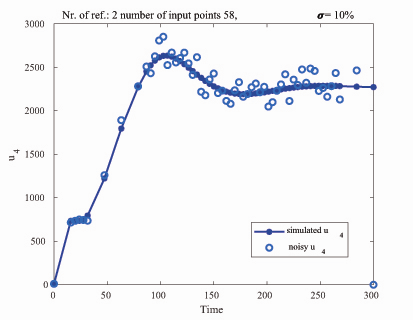}}  &
  {\includegraphics[scale=0.3, clip=]{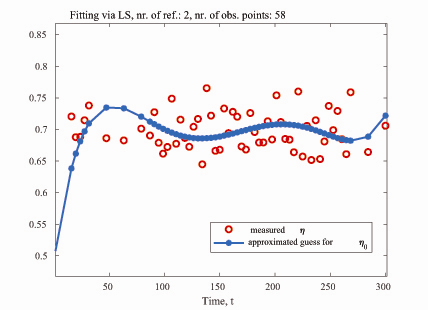}} &
  {\includegraphics[scale=0.3, clip=]{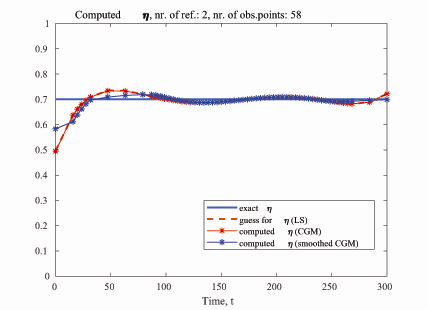}} \\
 ${u_4}_\tau$  & LS fitting to $\eta_\tau(t)$ & $\eta_\tau(t)$  \\
 \hline
\end{tabular}
\end{center}
\caption{
  \emph{  Test 2. Left figures: simulated ${u_4}_\tau$ vs.  noisy ${u_4}_\tau$ on different adaptively refined time meshes. Here, noisy observed  data are presented by circles.   Middle figures: least squares fitting to noisy data for $\eta_\tau$. Right figures: results of ACGA on adaptively refined meshes.  Computations are done for  noise level $\sigma=10\%$ in $u_4$  and for $T_1 = 100$.}}
\label{fig:FIG11}
\end{figure}



\begin{figure}[h!]
\begin{center}
  \begin{tabular}{ccc}
    \hline
    \multicolumn{3}{c}{  $nr.ref.=0$ }\\
  {\includegraphics[scale=0.3, clip=]{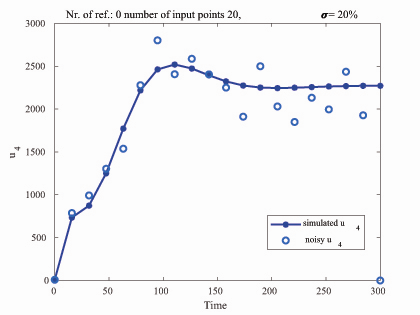}}  &
  {\includegraphics[scale=0.3, clip=]{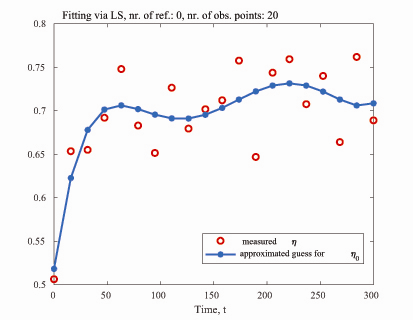}} &
  {\includegraphics[scale=0.3, clip=]{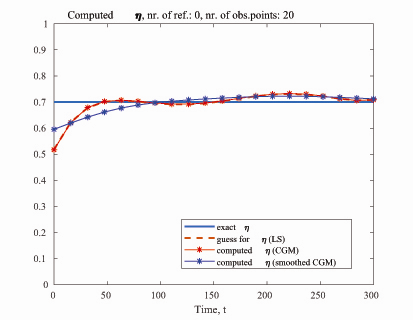}} \\
  ${u_4}_\tau$   & LS fitting to $\eta_\tau(t)$ & $\eta_\tau(t)$  \\
   \hline
  \multicolumn{3}{c}{  $nr.ref.=1$ }\\
  {\includegraphics[scale=0.3, clip=]{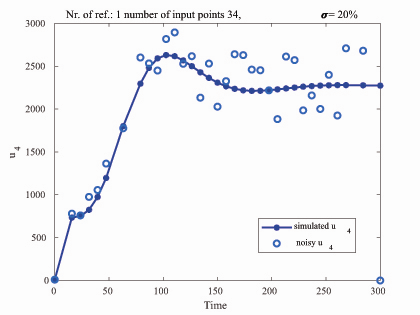}}  &
  {\includegraphics[scale=0.3, clip=]{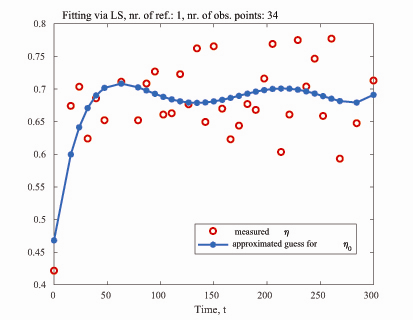}} &
  {\includegraphics[scale=0.3, clip=]{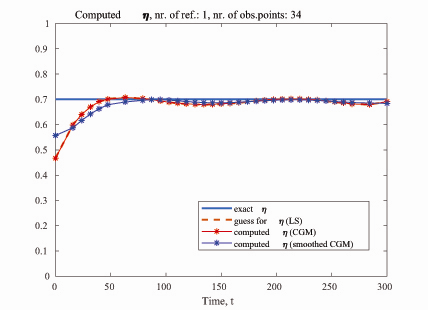}} \\
 ${u_4}_\tau$  & LS fitting to $\eta_\tau(t)$ & $\eta_\tau(t)$  \\
 \hline
   \multicolumn{3}{c}{  $nr.ref.=2$ }\\
  {\includegraphics[scale=0.3, clip=]{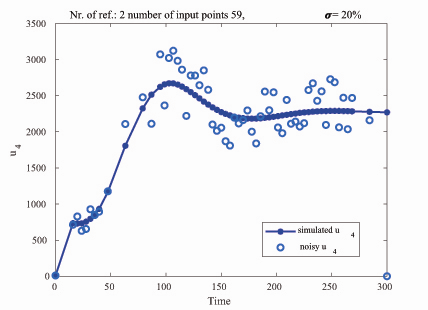}}  &
  {\includegraphics[scale=0.3, clip=]{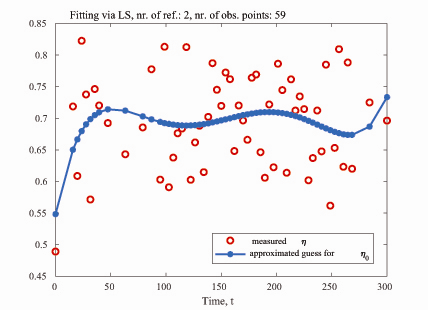}} &
  {\includegraphics[scale=0.3, clip=]{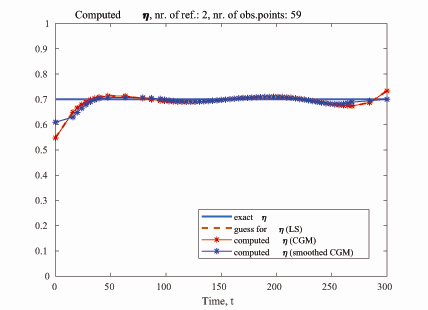}} \\
 ${u_4}_\tau$  & LS fitting to $\eta_\tau(t)$ & $\eta_\tau(t)$  \\
 \hline
\end{tabular}
\end{center}
\caption{
  \emph{  Test 2. Left figures: simulated ${u_4}_\tau$ vs.  noisy ${u_4}_\tau$ on different adaptively refined time meshes. Here, noisy observed  data are presented by circles.   Middle figures: least squares fitting to noisy data for $\eta_\tau$. Right figures: results of ACGA on adaptively refined meshes.  Computations are done for  noise level $\sigma=20\%$ in $u_4$  and for $T_1 = 100$.}}
\label{fig:test2b}
\end{figure}



\begin{figure}[h!]
\begin{center}
  \begin{tabular}{ccc}
    \hline
    \multicolumn{3}{c}{  $nr.ref.=0$ }\\
  {\includegraphics[scale=0.3, clip=]{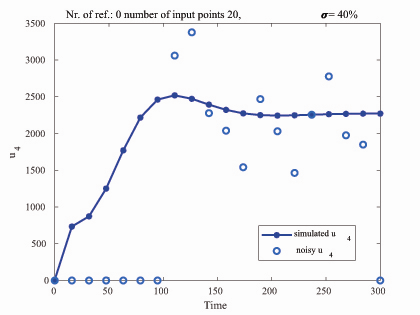}}  &
  {\includegraphics[scale=0.3, clip=]{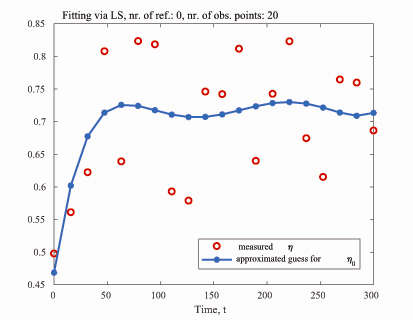}} &
  {\includegraphics[scale=0.3, clip=]{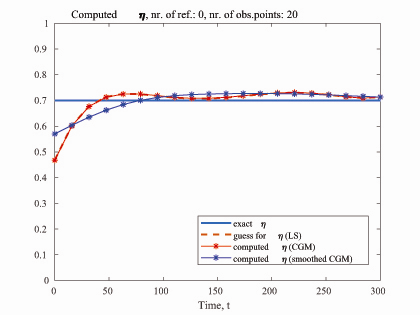}} \\
  ${u_4}_\tau$   & LS fitting to $\eta_\tau(t)$ & $\eta_\tau(t)$  \\
   \hline
  \multicolumn{3}{c}{  $nr.ref.=1$ }\\
  {\includegraphics[scale=0.3, clip=]{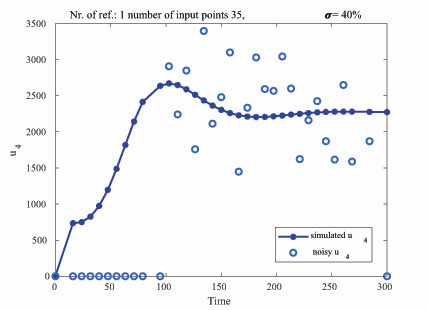}}  &
  {\includegraphics[scale=0.3, clip=]{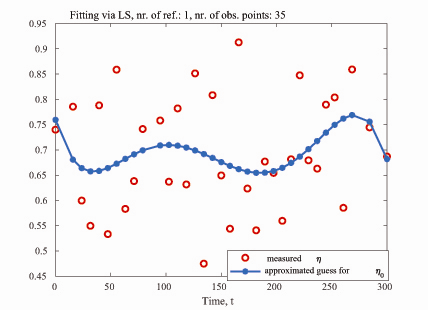}} &
  {\includegraphics[scale=0.3, clip=]{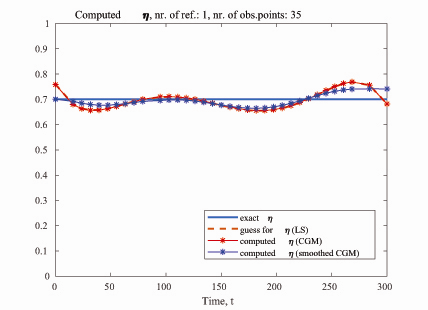}} \\
 ${u_4}_\tau$  & LS fitting to $\eta_\tau(t)$ & $\eta_\tau(t)$  \\
 \hline
   \multicolumn{3}{c}{  $nr.ref.=2$ }\\
  {\includegraphics[scale=0.3, clip=]{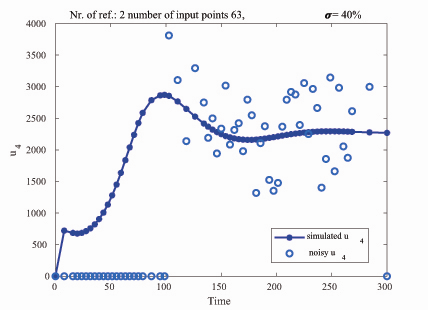}}  &
  {\includegraphics[scale=0.3, clip=]{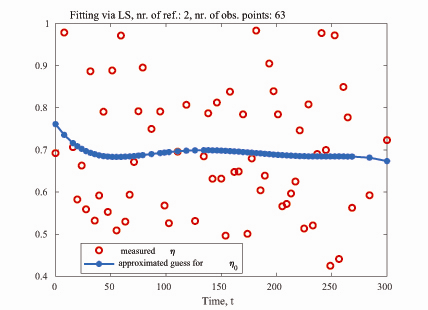}} &
  {\includegraphics[scale=0.3, clip=]{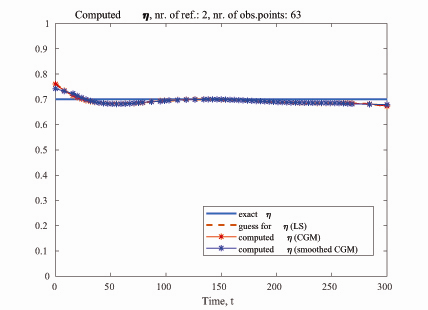}} \\
 ${u_4}_\tau$  & LS fitting to $\eta_\tau(t)$ & $\eta_\tau(t)$  \\
 \hline
\end{tabular}
\end{center}
\caption{
  \emph{  Test 2. Left figures: simulated ${u_4}_\tau$ vs.  noisy ${u_4}_\tau$ on different adaptively refined time meshes. Here, noisy observed  data are presented by circles.   Middle figures: least squares fitting to noisy data for $\eta_\tau$. Right figures: results of ACGA on adaptively refined meshes.  Computations are done for  noise level $\sigma=40\%$ in $u_4$  and for $T_1 = 100$.}}
\label{fig:FIG13}
\end{figure}


In this test we present numerical results of reconstruction of the model function $\eta(t)= 0.7$ from noisy observations of the function  $u_4(t)$ at the observation interval $[T_1, T_2]$. We again took $T_1= 25, 50, 100$, but number of observation points were $20$ at the time interval $[T_1, T_2] = [T_1, 300]$. We generate initial time partition $\mathcal J_{\tau}$ with equidistant time step $\tau = 300/19$.  Results of reconstruction of the model function $\eta(t)= 0.7$ for noise levels $\sigma = 5\%, 10\%, 20\%, 40\%$ in data $u_4(t)$ are presented in Table 2.  Figures \ref{fig:FIG6}-\ref{fig:FIG9} and \ref{fig:FIG10}-\ref{fig:FIG13} show results of reconstruction of the function $\eta(t)= 0.7$ for noise levels $\sigma = 5\%, 10\%, 20\%, 40\%$ in data $u_4(t)$ for $T_1 = 50$ and $T_1 = 100$, respectively.

 We again observe from the results of Table 2 and Figures \ref{fig:FIG6}-\ref{fig:FIG13}  that with local time-mesh refinements the reconstruction of the drug efficacy $\eta_\tau$ is  significantly improved compared to the reconstruction of $\eta_\tau$ obtained on initial non-refined time-mesh.

 \section{Conclusion}

 The time-adaptive optimization method for determination of drug
 efficacy in the mathematical model of HIV infection is presented.
 More precisely, first the time-dependent drug efficacy is determined
 at known coarse time partition using several known values of observed
 functions (usually 15-20 observations).  Then the time-mesh is
 locally refined at points where the residual $|R(\eta_\tau)|$ attains
 its maximal values and the drug efficacy is computed on a new refined
 time-mesh until the error in the reconstructed parameter $\eta$ is
 reduced to the desired accuracy.  Numerical experiments show
 efficiency and reliability of proposed adaptive method on
 reconstruction of different model functions $\eta$ from noisy
 observed virus population function.

The proposed new time-adaptive method can eventually be used by clinicians to
determine the drug-response for each treated individual. The exact
knowledge of the personal drug efficacy can aid in the determination
of the most suitable drug as well as the most optimal dose for each
person, in the long run resulting in a personalized treatment with
maximum efficacy and minimum adverse drug reactions.


\section*{Acknowledgment}

The research of the first author is supported by the Swedish
Research Council grant VR 2018-03661. The research of the second
author was supported by the  Russian Foundation for Basic Research
(grant 17-01-00636) and by the project N 0314-2018-0011.


\bigskip
\begin{thebibliography}{99}


\bibitem{BJ} W.~Bangerth and A.~Joshi, Adaptive finite element methods for the solution of inverse problems in optical tomography, \emph{Inverse Problems}, 24, 034011, 2008.

\bibitem{BR} R.~Becker and R.~Rannacher, An optimal control approach to \emph{a posteriori} error estimation in finite element method, \emph{Acta Numerica}, 10, pp.1--102, 2001.

\bibitem{Beilina1} L.~Beilina, Adaptive finite element method for a coefficient inverse problem for the {M}axwell's system, \emph{Applicable Analysis}, 90, pp.1461--1479, 2011.


\bibitem{Sprg-1} L.~Beilina and I.~Gainova, Time-adaptive FEM  for distributed parameter identification in biological models,  \emph{Applied Inverse Problems}, Springer Proceedings in Mathematics \& Statistics, 48, pp.37--50, 2013.

\bibitem{Sprg-2} L.~Beilina and I.~Gainova, Time-adaptive FEM  for distributed parameter identification  in mathematical model of HIV infection with drug therapy,  \emph{Inverse Problems  and Applications}, Springer Proceedings in Mathematics \& Statistics, 120, pp.111--124, 2015.

\bibitem{BJ2} L.~Beilina and C.~Johnson, \emph{A posteriori} error estimation in computational inverse scattering, \emph{Mathematical Models and Methods in Applied Sciences}, 15, pp.23--37, 2005.

\bibitem{BookBK} L.~Beilina, M~.V.~Klibanov, \emph{Approximate global convergence and adaptivity for Coefficient Inverse Problems}, Springer, New York, 2012.

\bibitem{BKK} L.~Beilina, M.~V.~Klibanov and M.~Yu.~Kokurin, Adaptivity with relaxation for ill-posed problems and global convergence for a coefficient inverse problem, \emph{Journal of Mathematical Sciences}, 167, pp.279--325, 2010.

\bibitem{Boch2012} G.~Bocharov, V.~Chereshnev, I.~Gainova, S.~Bazhan,
B.~Bachmetyev, J.~Argilaguet, J.~Martinez and A.~Meyerhans, Human
Immunodeficiency Virus Infection: from Biological Observations to
Mechanistic Mathematical Modelling, \emph{Mathematical Modelling of
Natural Phenomena}, 7(5), pp.78--104, 2012.

\bibitem{Burden} Richard~L.~Burden, J.~Douglas Fair\'es, Numerical Analysis, 9th Edition, Brooks/Cole

\bibitem{Cher2013}  V.~A.~Chereshnev, G.~A.~Bocharov, S.~I.~Bazhan,
B.~Bachmetyev, I.~A.~Gainova, V.~A.~Likhoshvai, J.~M.~Argilaguet,
J.~.P.~Martinez, J.~A.Rump, B.~Mothe, C.~Brander and A.~Meyerhans,
Pathogenesis and Treatment of HIV Infection: The Cellular, the
Immune System and the Neuroendocrine Systems Perspective,
\emph{International Reviews of Immunology}, 32(3), pp.282--306,
2013.
\\ http://informahealthcare.com/doi/abs/10.3109/08830185.2013.779375





\bibitem{Johnson} K.~Eriksson, D.~Estep, P.~Hansbo, C.~Johnson, \emph{Computational differential equations},
Cambridge University Press, 1996.

\bibitem{Martin} M. Eriksson, Parameter identification in a mathematical model of HIV infection with drug therapy, \emph{Master's thesis},  http://hdl.handle.net/2077/54664






 \bibitem{KB} N.~Koshev and L.~Beilina,  An adaptive finite element method for Fredholm integral equations of the first kind and its verification on experimental data,  in the Topical Issue ”Numerical Methods for Large Scale Scientific Computing” of CEJM, 11(8), 1489-1509,  2013.


\bibitem{Patrick} G.~L.~Patrick, \emph{An introduction to medicinal chemistry}, Fifth Ed., Oxford
University Press, Oxford, 2013.


\bibitem{Sriv09} P.~K.~Srivastava, M.~Banerjee, and P.~Chandra, Modeling the drug therapy for HIV infection, \emph{Journal of Biological Systems}, 17(2), pp.213--223, 2009.




\end{thebibliography}
\end{document}